\documentclass[11pt]{article} 
 
\usepackage{amsmath} 
\usepackage{amsthm} 
\usepackage{amssymb} 
\usepackage{amsfonts} 
 
\input xypic

\newtheorem{The}{Theorem}[section] 
 
\newtheorem{Pro}[The]{Proposition} 
\newtheorem{Exa}[The]{Example} 
\newtheorem{Rem}[The]{Remark} 
\newtheorem{Lem}[The]{Lemma} 
\newtheorem{Def}[The]{Definition}

\def\proof{\vspace{2ex}\noindent{\bf Proof.} }

\def\endproof{\relax\ifmmode\expandafter\endproofmath\else 
\unskip\nobreak\hfil\penalty50\hskip.75em\hbox{}\nobreak\hfil\bull 
{\parfillskip=0pt \finalhyphendemerits=0 \bigbreak}\fi} 
\def\endproofmath$${\eqno\bull$$\bigbreak} 
\def\bull{\vbox{\hrule\hbox{\vrule\kern3pt\vbox{\kern6pt}\kern3pt\vrule}\hrule}} 
 
\def\dirac{\partial}

\newcommand{\C}{{\mathbb C}} 
\newcommand{\R}{{\mathbb R}} 
\newcommand{\Z}{{\mathbb Z}} 
\newcommand{\Q}{{\mathbb Q}} 
\newcommand{\AAA}{{\mathbb A}} 
\newcommand{\A}{{\cal A}} 
\newcommand{\E}{{\cal E}} 
\newcommand{\G}{{\cal G}} 
\newcommand{\M}{{\cal M}} 
\newcommand{\B}{{\cal B}} 
 
\newcommand{\ind}{\mathfrak{i}}

\newcommand{\la}{\langle} 
\newcommand{\ra}{\rangle} 
\newcommand{\ba}{\begin{eqnarray}} 
\newcommand{\na}{\end{eqnarray}} 
\newcommand{\beq}{\begin{equation}} 
\newcommand{\eeq}{\end{equation}}

\newcommand{\s}{\mathfrak{s}} 
\renewcommand{\t}{\mathfrak{t}} 
\newcommand{\spinc}{\mathrm{Spin}^c}

\title{Exact triangles in Seiberg-Witten-Floer theory. Part IV: $\Z$-graded 
monopole  homology}  
\author{Matilde Marcolli and Bai-Ling Wang} 
\date{}
 
\begin{document} 
\maketitle

\section{Introduction} 
 
The content of this paper is part of our ongoing project  
of establishing the surgery formulae for 
Seiberg-Witten-Floer theory. As explained in the introductory part of  
\cite{CMW}, if $K$ is a knot in a homology 3-sphere 
$Y$,  we expect the Floer homologies of $Y$ and of the manifolds $Y_1$  
and $Y_0$ obtained by Dehn surgery on $K$ with framing $1$ and $0$, 
respectively, to be related by an exact triangle 
$$ 
\diagram 
HF^{SW}_{*} (Y_1, g_1) \rrto^{I_{*}}  
&  & HF^{SW}_{*}(Y, g) 
\dlto_{L_{*}} \\  
& \bigoplus_k HF^{SW}_{(*)}(Y_0, \s_k) \ulto_{\Delta_{(*)}}  
& \\ 
\enddiagram 
$$ 
 
There are several subtleties involved in the precise definition of the term 
$$\bigoplus_k HF^{SW}_{(*)}(Y_0, \s_k)$$ in the previous diagram. In this 
paper we explain the precise meaning of this term.  
 
More specifically, we analyze some aspects of the construction 
of the Seiberg-Witten-Floer homology of the manifold $Y_0$, obtained 
by 0-surgery on the knot  
$K$ in the homology sphere $Y$. The manifold $Y_0$ has the homology of 
$S^1\times S^2$. For manifolds with $b_1(Y_0)>0$ we know (cf. \cite{CW}, 
\cite{Ma}, \cite{MW}, \cite{RW}) that there is a well defined $\Z_\ell$ graded 
Seiberg-Witten-Floer homology $HF^{SW}_*(Y_0, \s)$, for every choice 
of a non-trivial $\spinc$ structure $\s \in {\cal S}(Y_0)$, 
where the integer $\ell=\ell(\s)$ satisfies 
$$ c_1(L) ( H_2(Y_0,\Z))  = \ell \Z, $$
and ${\cal S}(Y_0)$ denotes the set of $\spinc$ structures on $Y_0$. 
Here $L=\det W$ is the determinant line bundle of the spinor bundle 
$W$ associated to the $\spinc$ structure $\s$. We use the equivalent 
notation $L=\det \s$. Here we are assuming that $c_1(L)$ is
non-torsion, hence $\ell\neq 0$. 
 
In this paper we shall discuss several issues regarding the Floer homology  
$HF^{SW}_*(Y_0, \s)$. The first is connected to the integer lift of 
the $\Z_\ell$ graded Floer homology.  
Analogous constructions of integer lifts of Floer 
homologies were derived by Fintushel and Stern \cite{FS}, in the case 
of instanton homology, and by Weiping Li \cite{Li}, in the case of 
symplectic Floer homology. We follow  
closely the construction of \cite{FS} 
and show that, in our case, there is a well defined integer  
lift of $HF^{SW}_{*, (\omega) }(Y_0, \s)$ of the $\Z_\ell$-graded Floer 
homology, here $\omega\in \R$ is a regular value of the 
Chern-Simons-Dirac functional on the infinite cyclic cover space   
of the gauge equivalence classes of connections and spinor sections. By 
 studying the Chern-Simons-Dirac function on this infinite cyclic cover space, 
we will define an integer lift $\ind_{Y_0}^{(\omega)}$ 
of the indices of the critical points. We thus form 
a chain complex $C_{*}^{(\omega)}(Y_0, \s)$ depending on  
$\omega \in \R$.  
For any $n \in 
\Z_{\ell}$, the original $\Z_{\ell}$-graded Seiberg--Witten--Floer 
chain complex satisfies 
$C_n(Y_0, \s) = \bigoplus_{k\in \Z} C_{j +k\ell}^{(\omega)} (Y_0, \s)$  
where $j=n (mod \ell)$.  
In general, after defining a suitable boundary operator on
$C_{*}^{(\omega)}(Y_0, \s)$, we observe that the resulting homology groups 
$HF_{*, (\omega)}^{SW} (Y_0, \s)$ do not satisfy the simple relation
$\bigoplus_{k\in \Z} HF_{* + k\ell, (\omega)}^{SW}(Y_0, \s) =
HF_*^{SW}(Y_0, \s)$. However, the Floer homologies 
$HF_{* + k\ell, (\omega)}^{SW}(Y_0, \s)$ and $HF_*^{SW}(Y_0, \s)$ are
related via a spectral sequence determined by a  
filtration of the chain complex $C_{*}(Y_0, \s)$. This spectral
sequence converges to $HF^{SW}_*(Y_0, \s)$, and the $E^1$ term 
coincides with $HF^{SW}_{*, (\omega)}(Y_0, \s)$. 
For instance, in the particular case where all the higher
differentials in the spectral sequence   
are zero, then we simply have,  
for $m \in \{ 0, 1, \cdots, \ell -1\}$, 
\[ 
HF_m^{SW}(Y_0, \s) = 
\bigoplus_{n\in \Z} HF^{SW}_{m+ n\ell, (\omega)}(Y_0, \s). 
\] 
The arguments we develop in this section extend easily to the general
case of any 3-manifold $(Y_0, \s)$ with non-trivial rational homology
and with a $\spinc$-structure satisfying $c_1(\s) (H_2(Y_0, \Z)) =
\ell \Z \neq 0$.  
 
The second main issue discussed in this paper is a different
construction of a $\Z$-graded Seiberg--Witten--Floer homology, where the 
integer lift of the original $\Z_{\ell}$-graded Seiberg--Witten--Floer 
homology is determined by the exact triangle, as derived in \cite{CMW} 
\cite{MW2} \cite{MW3}. 
Under the surgery identifications of the Seiberg--Witten monopoles 
on $Y$ and $Y_0$, proved in \cite{CMW}, there exists a one-to-one map 
\ba 
 j: \qquad \bigsqcup_{\s} \M_{Y_0} (\s) \longrightarrow  \M_{Y, \mu}
\backslash  \M_{Y_1}. 
\na 
 Then the $\Z$-valued indices on $\M_{Y, \mu}$  
define a $\Z$-lifting $\ind_{Y_0}^{(Y)}$ of the grading $\M_{Y_0}(\s)$. 
The corresponding chain complex $C_{(n)}(Y_0, \s)$ is 
defined to be $\Z\{a \in  \M_{Y_0} (\s) : \ind_{Y_0}^{(Y)}(a) =n\}.$ 
We will show that the restriction of the boundary operator 
on $C_*(Y_0, \s)$ to $C_{(n)}(Y_0, \s)$ is a well-defined boundary
operator, and the resulting homology groups $HF_{(*)}(Y_0, \s)$ satisfy 
\ba 
HF_n^{SW}(Y_0, \s) = 
\bigoplus_{k\in \Z} HF^{SW}_{(n+ k\ell)}(Y_0, \s). 
\label{oplus} 
\na 
The analysis of the splitting and gluing of the moduli spaces of flow lines 
in \cite{MW2} play an essential role in the proof of the  
direct sum formula (\ref{oplus}). 
It is precisely these $\Z$-graded monopole homology groups that
we use in the exact triangles \cite{MW3}.

The third main issue that we discuss in this paper is the construction of 
the Seiberg-Witten-Floer homology on a 3-manifold $Y_0$ with
$b_1(Y_0)>0$, in the case of a $\spinc$  
structure $\s_0$ with torsion $c_1(\det \s_0)$. This case was not included
in our original paper \cite{MW}. We show that, in this case, 
the boundary operator of the Floer complex is well defined only after 
separating the components of uniform energy in the moduli space ${\cal  
M}(a,b)$. This is done by introducing a Novikov complex with 
coefficients in $\Z[[t]]$, which depends on the  
non-trivial cohomology class of the perturbation for the Chern-Simons-Dirac 
functional. The component of this Novikov-Floer 
homology $HF_*(Y_0,\s_0,\Z[[t]])$ that we are interested in, for the
purpose of the exact triangle, is the evaluation 
$$ HF_{(*)}(Y_0,\s_0)=HF_*(Y_0,\s_0,\Z[[t]])|_{t=0}. $$

We also generalize the construction of 
$HF_*(Y_0,\s_0,\Z[[t]])$  
to a more natural setting over the coefficient field of formal  
Laurent series $\Q((t)) = \Q[[t]][t^{-1}]$.  
One reason for this 
definition is that, in this setting, we can associate a relative
Seiberg--Witten invariant for any 4-manifold $(X, \s)$ with boundary
$(Y_0, \s_{0})$, where $c_1(\s|_{Y_0})$ is a torsion class.  Here
we introduce a non-trivial 1-cycle $\Gamma$ in $Y_0$, 
which depends on the  non-trivial cohomology class of the perturbation.
The corresponding Seiberg-Witten-Floer homology is
denoted by $HF^{SW}_{\Gamma, *}(Y_0, \s_0, \Q((t)))$. 

Suppose given 4-manifold $(X, \s)$ with a cylindrical end modeled on
 $(Y_0, \s_0 = \s|_{Y_0})$, where $c_1(\s_0)$ is a torsion cohomology
class in $H^2(Y_0. \Z)$. We study
the perturbed 4-dimensional Seiberg-Witten  equations
on $(X, \s)$, the perturbation for $(X, \s)$
is compactible with the perturbation on $(Y_0, \s_0)$
used to the construction of $HF^{SW}_{\Gamma, *}(Y_0, \s_0, \Q((t)))$.
For simplicity, we assume that $\Gamma$ bounds a relative
2-cycle in $X$, then 
the relative Seiberg--Witten 
invariant $SW_X(\s)$  
takes values in $HF^{SW}_{\Gamma, *}\bigl(Y_0,\s_0,\Q((t))\bigr)$. We derive 
the gluing formulae for this case. 
When $\Gamma$ does not bound  any relative
2-cycle in $X$, the relative invariant takes values   in
the Seiberg-Witten-Floer homology $(Y_0, \s_0)$
over the coefficient field of formal  
Laurent series on $H^1(Y_0, \Z)$.  The gluing formulae along the torsion $\spinc$  structure can be formulated similarly.

  As an example, our gluing formula
can be applied to show that, for 
certain particular choice of metric and perturbation on the three torus 
$S^1\times T^2$ with the trivial $\spinc$ structure $\s_0$, we have 
\[ 
HF_*\bigl(S^1 \times T^2, \s_0, \Q((t))\bigr) = \Q((t)), \]
and for $D^2\times T^2$ with the cylindrical end
modeled on $S^1\times T^2$, we get
\[
  SW_{D^2 \times T^2} (\s_0) = \displaystyle{\frac 
{1}{t  - t^{-1}}}. 
\] 
In the more general setting, these gluing formulae provide useful tools 
in the study of the structure of the Seiberg--Witten--Floer  
homology groups (cf. \cite{vicente}). The relative invariant
$SW_{D^2 \times T^2} (\s_0)$ can also be applied to knot surgery of
4-manifolds to give some of  the gluing formulae obtained 
in \cite{FS2} and \cite{MMS}.

\vskip .2in 
 
{\bf Acknowledgments} Both the authors benefited of conversations 
with Ron Wang on some of the issues discussed in this paper. 
We thank Vicente Mu\~noz
for his help in formulating the gluing formulae. We also thank the
referee for useful comments. We thank the Max Planck Institut f\"ur
Mathematik for the kind hospitality and for support. The first author
is partially supported by NSF grant DMS-9802480. The second author is
partially supported by ARC Fellowship.

\section{ Seiberg-Witten-Floer homology }

In the following two sections we present two different constructions of 
an integer lift of the Seiberg--Witten--Floer homology. The first case
applies to the general setting where  $Y_0$ is a
compact 3-manifold with $b_1(Y_0)>0$, and the $\spinc$-structure is 
non-torsion. The second applies to the special case where $Y_0$ is 
obtained by zero-surgery on a knot in a homology 3-sphere $Y$, again
with a non-torsion $\spinc$-structure.
In this particular case, the integer lift of the Seiberg--Witten--Floer
homology that we consider appears naturally in the context of the 
exact triangle formula. 

We give some preliminary definitions.
 
Let ${\cal M}_{Y_0}(\s)$ be the moduli space of gauge classes of 
solutions of suitably  
perturbed Seiberg--Witten equations on a 3-manifold $(Y_0, \s)$, 
with $b_1(Y_0) >0$ and $c_1(\s)$ a non-torsion element. A 
suitable class of perturbations that achieves transversality has been 
introduced in \cite{CMW}. The 
perturbation can be chosen so that all the solutions in $\M_{Y_0}(\s)$  
are irreducible critical points. Thus, under the choice of a generic 
perturbation, $\M_{Y_0}(\s)$ is a compact, oriented, 
0-dimensional manifold, cut out transversely by the equations inside the 
configuration space $\B =\A / \G$, where $\A$ is the 
space of pairs formed by a $U(1)$-connection on $det(\s)$ and a spinor
section of $W$, and $\G$ is the group of gauge transformations 
on $(Y_0, \s)$.  
 
Recall that we have the Chern--Simons--Dirac functional on $\A$,
defined as  
\[\begin{array}{lll} 
 CSD(A,\psi)& =& - \displaystyle{\frac{1}{2}\int_{Y_0} }(A-A_0)\wedge (F_A +  
F_{A_0}-2\sqrt{-1} *\rho) \\[2mm] 
&& + \displaystyle{\int_{Y_0} }\la \psi,\partial_A\psi \ra dvol_{Y_0},  
\end{array} 
\] 
with $\rho$ a co-closed 1-form on $Y_0$. In order to achieve 
transversality in the moduli spaces of flow lines, we consider an  
additional perturbation of $CSD$ by functions 
$ U(\tau_1,\ldots,\tau_N)+ V(\zeta_1,\ldots,\zeta_K) $ 
as explained in \cite{CMW}. Here $\tau_j$ and $\zeta_j$ 
are functions on $\A$ respectively associated to a complete basis
$\{\mu_j\}_{j=1}^\infty$  
for the co-closed $L^2$ one-forms on $Y_0$ and a complete 
basis $\{\nu_j\}_{j=1}^\infty$ for the imaginary-valued one-forms on $Y_0$, 
\[\begin{array}{l} 
 \tau_j (A, \psi) = \tau_j (A-A_0) =\displaystyle{ \int_{Y_0}} (A-A_0)
\wedge * \sqrt{-1}\mu_j\\[2mm]  
\zeta_j (A, \psi) = \zeta_j (\psi, \psi)= 
\displaystyle{ \int_{Y_0} }\la \nu_j. \psi, \psi \ra dvol_{Y_0}. 
\end{array} 
\] 
The pair $(U, V)$ is chosen from a subspace of functions in
\[
\bigcup_{N\ge b_1, K>0} C^\infty(\R^N, \R) 
\times C^\infty (\R^K, \R)\]
completed to 
a Banach space in the Floer $\epsilon$-norm (as specified in \cite{CMW}).
This perturbation term is gauge invariant by construction, and $CSD$ 
changes by 
$$  
CSD(\lambda.(A, \psi))-CSD(A,\psi) 
= \la (8\pi^2   c_1(L) +4\pi [*\rho] ) \cup [\lambda],[Y_0]\ra, 
$$  
where $[\lambda]\in H^1(Y_0,\Z)$ determines the connected component of  
${\cal G}$ that contains $\lambda$, and is represented by the closed 1-form 
$ \frac{1 }{2\pi \sqrt{-1}} \lambda^{-1}d\lambda.$ 

With this choice of perturbation, the Seiberg--Witten equations are of
the form
\[
\left\{ \begin{array}{l}
*F_A = \sigma (\psi, \psi) +\displaystyle{
 \sum_{j=1}^{N} \frac{\partial U}{\partial \tau_j}}\mu_j \\[2mm]
\dirac_A (\psi)+ \displaystyle{\sum_{j=1}^{K} \frac{\partial
V}{\partial \zeta_j}} 
\nu_j.\psi =0,
\end{array}\right.
\]
for the critical points of the Chern--Simons--Dirac functional, and
\[ 
\left\{ \begin{array}{lll}
\displaystyle{\frac{\partial A}{\partial t} }&=& -*F_A + \sigma (\psi, \psi) + 
\sum_{j=1}^{N} \displaystyle{\frac{\partial U}{\partial \tau_j}} \mu_j \\[2mm]
\displaystyle{\frac{\partial \psi}{\partial t}}& =& -\dirac_A \psi 
- \sum_{j=1}^{K} \displaystyle{\frac{\partial V}{\partial \zeta_j}}\nu_j.\psi
\end{array}\right.
\]
for the corresponding downward gradient flow.

The Seiberg--Witten--Floer homology groups 
$HF^{SW}_*(Y_0, \s)$ (cf. \cite{CW}, \cite{Ma}, \cite{MW})  
are the homology groups 
of $(C_*(Y_0, \s), \partial)$, where  
$C_*(Y_0, \s)$ is generated by the critical points of the perturbed
$CSD$ on $\B$. The entries of the boundary operator $\partial$ are defined 
by counting the points in the zero dimensional components of the moduli 
space of unparameterized flow lines of the perturbed functional $CSD$ on
$\B$. A considerable amount of technical work goes into checking that
the moduli spaces of flow lines have all the desired properties that
make this definition rigorous, and we refer the reader to \cite{MW}
for a precise account. The resulting Floer homology $HF^{SW}_*(Y_0, \s)$  
is $\Z_\ell$-graded with $\ell \Z = c_1(\s) (H_2(Y_0, \Z))$. This is due to 
an ambiguity in the index formula computed in terms of the spectral 
flow of the Hessian operator for $CSD$ around a loop in $\B$. To
understand this ambiguity, let  
$(A, \psi)$ represent a critical point $a$ of $CSD$ on $\B$, let $\lambda$ 
be a gauge transformation whose class $[\lambda] \in H^1(Y_0, \Z)$ 
is non-trivial. Then the spectral flow of the Hessian operator $H$ from 
$\lambda. (A, \psi)$ to $ (A, \psi)$ can be calculated from the 
index formula for the linearization of the 4-dimensional
Seiberg--Witten equations on $Y_0 \times S^1$. This gives
\ba 
SF(H)_{\lambda. (A, \psi)}^{(A, \psi)} = \la [\lambda] \cup c_1(\s), [Y_0]\ra 
\in \ell \Z. 
\label{index:ell} 
\na

\begin{Rem} 
Notice that the periodicity $\ell$ is an even number 
(cf. \cite{CMW}). In fact, the map 
${\cal S}(Y_0)\to H^2(Y_0,\Z)$, given by 
$\s \mapsto c_1(\det \s)$, 
is equivariant with respect to the action of $H^2(Y_0,\Z)$, so that we 
have $ \det (W\otimes H) = L\otimes H^2$, 
with $L=\det(W)=\det(\s)$ and $H\in H^2(Y_0,\Z)$, cf. \cite{Mabook}.  
\end{Rem}

There is a natural and non-trivial way to lift  
these $\Z_\ell$-graded homology groups 
to $\Z$-graded homology groups, so that the  Euler characteristic number 
agrees with the original one. Following the idea of  
\cite{FS},  we will discuss  this integer lift and 
construct a spectral sequence which has the integer lift as  
$E^1$ term  and converges to the original $\Z_\ell$-graded Floer
homology.  We shall then  
consider a different lifting of the $\Z_\ell$-graded homology 
groups which is relevant to the exact triangles considered in 
\cite{CMW}, \cite{MW3}, and \cite{MW4}. 
 
\section{$\Z$-graded homology groups and the spectral sequence} 

The construction we present in this section holds in the general case
of a 3-manifold $Y_0$ with $b_1(Y_0)>0$ and a $\spinc$-structure with
non-torsion class $c_1(L)$.
 
As explained in \cite{MW}, there is a cyclic covering of $\B$, 
obtained by taking the quotient of $\A$, the space of 
$U(1)$-connections and spinors, with respect to the subgroup $\G_\ell$  
of the gauge group $\G$ given by 
\[ 
 \G_\ell = \{ \lambda \in \G \ | \ \la c_1(L)\cup [\lambda], [Y_0] 
\ra=0 \},  
\] 
where $\ell$ satisfies $\ell \Z = c_1(L) (H_2(Y_0,\Z) )$.  
This subgroup depends on $c_1(L)$.
 
The space $\B$ has the homotopy type of   
$\C P^\infty \times K(H^1(Y_0,\Z), 1)$,  
hence it has a universal covering obtained by 
taking the quotient of $\A$ by the identity component of the gauge 
group. The resulting space $\tilde \B$ covers $\B$ with fibers 
$H^1(Y_0,\Z)$.  
 Define $H_\ell$ to be  
\[ H_\ell=\{ h\in H^1(Y_0,\Z) | \la c_1(L)\cup h, [Y_0] \ra =0 \}. 
\] 
The group $H_\ell$ also depends on $c_1(L)$.
Then the  space $\B_\ell=\A /\G_\ell$ is a covering of $\B$ with fiber 
$ H^1(Y_0,\Z)/ H_\ell \cong \Z$. Hence $\B_\ell$ is an infinite  
cyclic covering space of $\B$.

The perturbed Chern--Simons--Dirac functional is a real valued 
functional $CSD: \B_\ell \to \R$ on the covering space 
$\B_\ell$.  
The critical manifold $\tilde \M_{Y_0}(\s)$ is a $\Z$-covering of 
$\M_{Y_0}(\s)$. The critical values form a discrete set in $\R$, which 
is a finite set mod $\Z$. Let 
$\Omega \subset \R$ denote the set of regular values. 
Let $\omega\in \Omega$ be a regular value. Given any point $a\in 
\M_{Y_0}(\s)$, there is a unique element $  a^\omega$ in the fiber 
$\pi^{-1}(a)$ in   
$\tilde \M_{Y_0}(\s)$ that satisfies  
$CSD( a^\omega)\in (\omega,\omega+8\pi^2 \ell).$ 
 
We have the following Lemma which 
shows that the relative indices on $\tilde \M_{Y_0}(\s)$, 
 defined by the spectral flow of the Hessian operator, 
take values in $\Z$. 
 
\begin{Lem}\begin{enumerate} 
\item 
Consider the Hessian operator $H_{A(t),\psi(t)}$ defined in 
\cite{CMW}. The 
spectral flow of the operator $H_{A(t),\psi(t)}$ around a loop in 
$ \B_\ell$ is zero, hence the relative index of any two points $\tilde a$  
and $\tilde b$ in $\tilde \M_{Y_0}(\s)$ is well defined in $\Z$. 
\item Let $\tilde a$ be a critical point in  $\tilde \M_{Y_0}(\s) \in
\B_\ell$,  
and let $\lambda \in \G/\G_\ell$ be a gauge transformation. We
have the following identities:  
\[\begin{array}{l} 
CSD(\lambda(\tilde a)) -CSD(\tilde a) = 8\pi^2\la 
   [\lambda ] \cup c_1(\s), [Y_0]\ra;\\[2mm] 
SF(H_{(A(t),\psi(t))})|_{\lambda(\tilde a)}^{\tilde a} =  \la 
   [\lambda ] \cup c_1(\s), [Y_0]\ra, 
\end{array} 
\] 
with $[\lambda] = [\frac{1 }{2\pi \sqrt{-1}} \lambda^{-1}d\lambda]$.
Here  $SF(H_{(A(t),\psi(t))})|_{\lambda(\tilde a)}^{\tilde a}$ denotes
the spectral 
flow of the Hessian operator $H_{(A(t),\psi(t))}$ along a path from
$\lambda(\tilde a)$ to $\tilde a$ in $\B_\ell$. 
\end{enumerate} 
\label{index:shift} 
\end{Lem} 
\proof 
By the Atiyah-Patodi-Singer index theorem, we have
\[ 
SF(H_{(A(t),\psi(t))})|_{\lambda(\tilde a)}^{\tilde a} = Index
(\displaystyle{\frac{\partial}{\partial t} +   
H_{(A(t),\psi(t))}}). 
\] 
This index calculates the virtual dimension of the 
4-dimensional Seiberg--Witten monopole moduli space on $Y_0\times S^1$, with 
the $\spinc$ structure obtained by gluing the $\spinc$ 
structure $\s$ on $Y_0\times \R$ along the two ends 
with the gauge transformation $\lambda$. By the index formula 
for the Seiberg--Witten monopoles, we obtain 
\[ 
\begin{array}{lll} 
&&Index (\displaystyle{\frac{\partial}{\partial t} +  
H_{(A(t),\psi(t))}})\\[2mm] 
&=&\displaystyle{ \frac{1}{2\pi\sqrt{-1}} \int_{Y_0} }c_1(L)\wedge
\lambda^{-1}d\lambda\\[3mm]  
&=& \la [\lambda ] \cup c_1(\s), [Y_0]\ra. 
\end{array} 
\] 
The remaining claims are direct consequence of this index formula. 
\endproof 
 
Upon fixing a base point $\tilde a_0$ in $\tilde \M_{Y_0}(\s)$, we can  
define the following $\Z$-lifting of the $\Z_\ell$-grading of the 
elements of $\M_{Y_0}(\s)$. 
 
\begin{Def} 
We define the grading of elements in $\M_{Y_0}(\s)$ as the relative 
index of the $\omega$-lifting in $\tilde \M_{Y_0}(\s)$, 
$$ \ind_{Y_0}^{(\omega)} (a)= \ind_{Y_0}(a^\omega) =
SF(H_{(A(t),\psi(t))})|_{a^\omega}^{\tilde a_0}. $$  
\end{Def} 
 
This definition of the grading depends on the choice of the base point  
$  \tilde a_0$ and on the choice of a regular value $\omega$. Notice
that we can reduce the choice of $\tilde a_0$ to  
the choice of a base point $a_0$ in $\M_{Y_0}(\s)$. 
  From now on, we fix a base point $a_0$  in $\M_{Y_0}(\s)$. 
Then $\tilde a_0$ is chosen to be  
the unique critical point in $\pi^{-1}(a_0)$ with  
$CSD(\tilde a_0) \in (\omega, \omega + 8 \pi^2)$. 
It is easy to see, from 
the definition, that  we have $\ind_{Y_0}^{(\omega)} =
\ind_{Y_0}^{(\omega')}$ whenever
$\omega$ and $\omega'$ are connected by a path in the set $\Omega$ of
regular values. 
Moreover, we have $\ind_{Y_0}^{(\omega +  8\pi^2 k\ell)} =
\ind_{Y_0}^{(\omega)}$, for any $k\in \Z$.

\begin{Def} 
For any $q \in \Z$, we define 
\[ 
 \M_{Y_0, q }^{(\omega)}(\s)=  
\{ a\in \M_{Y_0}(\s) \ | \ind_{Y_0}^{(\omega)}(a)=q \}.\] 
The $q$-chains in the $\Z$-graded Floer complex are the elements of 
the abelian group $C_q^{(\omega)}(Y_0,\s)$  
generated by the monopoles in $\M_{Y_0, q}^{(\omega)}(\s)$. 
The boundary operator 
$$ \partial_{q}^{(\omega)} : C_{q}^{(\omega)}(Y_0,\s) 
 \to C_{q-1}^{(\omega)} (Y_0,\s) $$ 
is defined as 
\[ 
\partial_q^{(\omega)}(a) = \sum_{b\in \M_{Y_0, q-1 }^{(\omega)}(\s)} 
\# \bigl( \hat \M^0 (a, b)\bigr) b, 
\] 
where $\hat \M^0 (a, b)$ is the zero dimensional components of the 
moduli space of  
unparameterized flow lines in $\hat\M(a, b)$. The compactness theorem
(cf. \cite{MW})
tells us that $\hat \M^0 (a, b)$ is an oriented, compact, 0-dimensional
manifold. Thus, the coefficient $\#(\hat \M^0 (a, b))$ is well-defined 
as the algebraic sum of points in $\hat \M^0 (a, b)$.  
\label{(q)} 
\end{Def} 
 
The following Lemma shows that we have $\partial_{q-1}^{(\omega)} 
\circ \partial_q^{(\omega)} =0$. The resulting homology groups 
are denoted as $HF^{SW}_{*, (\omega)}(Y_0, \s)$, $*\in \Z$.  
 
\begin{Lem} \begin{enumerate} 
\item For $n\in \Z_\ell$ and $q\in \Z$ with $q=n (mod \ \ell)$, we have  
\[ 
C_n (Y_0, \s)=  \bigoplus_{k\in \Z}  C_{q+ k \ell}^{(\omega)}(Y_0,\s). 
\] 
\item Under the decomposition as 
above, the boundary operator $\partial_n$ on 
$C_n(Y_0, \s)$ can be expressed as 
\[ 
\partial_n =  
 \left( \begin{array}{ccccccc} 
\partial_*^{(\omega)} &\partial_{*, 1}^{(\omega)} 
   & \partial_{*, 2}^{(\omega)}&\cdots & \cdots &\partial_{*,
*}^{(\omega)}  \\  
0 & \partial_{*+\ell} ^{(\omega)}&  \partial_{*+\ell, 1}^{(\omega)}
&\cdots &\cdots &\partial_{*+\ell, *}^{(\omega)} \\  
0 & 0 & \partial_{*+2\ell}^{(\omega)}& \cdots &\cdots &  
 \partial_{*+2\ell, *}^{(\omega)}\\ 
0&0&0& \ddots && \partial_{*+3\ell, *}^{(\omega)}\\ 
\vdots & \vdots &\vdots &&  \ddots & \vdots  \\ 
0&0&0 &\cdots &\cdots & \partial_{*+ *\ell}^{(\omega)} 
 \end{array} \right). 
\] 
The meaning of this  matrix notation can be
explained more precisely as follows. Assume that $a$ is a generator in
$C_q^{(\omega)}(Y_0, \s)$. Upon regarding $a$ as 
a generator of $C_n (Y_0, \s)$ for $n= q(mod\ \ \ell)$,
we obtain
\[ 
\partial_n (a) = \partial_q^{(\omega)} (a) + \sum_{k>0} 
\partial_{q, k}^{(\omega)} (a), 
\] 
with $\partial_{q, k}^{(\omega)}: \ 
C_q^{(\omega)} (Y_0, \s) \rightarrow C_{q-1+k \ell}^{(\omega)}$ for $k>0$. 
In particular, the relation $\partial_{n-1}\circ \partial_n = 0$ implies 
that $\partial_{n-1}^{(\omega)}\circ \partial_n^{(\omega)} = 0$ is
also satisfied. 
\end{enumerate} 
\label{boundary:operator} 
\end{Lem}
 
\proof The first statement about the decomposition of the 
chain complex is obvious from the definition. Now we study 
the boundary operator under this decomposition. 
 
For any $k<0$, $a\in \M_{Y_0, q}^{(\omega)}(\s)$, and  
$b\in \M_{Y_0, q-1 + k\ell}^{(\omega)}(\s)$, we shall prove that 
the 0-dimensional components $\hat \M^0(a, b)$ in 
$\hat \M (a, b)$ is empty, hence the entry of the boundary operator is 
trivial, $\la b, \partial_n(a) \ra =0$. 
 
Notice that we have 
$$ \ind_{Y_0}^{(\omega)} (a) = \ind_{Y_0} (a^\omega) = q \ \ \hbox{ and
} \ \ \ind_{Y_0}^{(\omega)} (b) = \ind_{Y_0} (b^\omega) = q-1 + k\ell. $$ 
Thus, the moduli space of flow lines on $\B_\ell$, 
$\hat\M (a^\omega, b^\omega)$ has virtual dimension $-kl>0$.  
There exists a unique element $[\lambda] \in \G/\G_\ell$, 
such that  
\[ 
\la [\lambda] \cup  c_1(L), [Y_0]\ra = -kl. 
\] 
This implies that we have $\ind_{Y_0} (\lambda (b^\omega)) =  
SF(H_{(A(t), \psi(t))})_{\lambda (b^\omega)}^{\tilde a_0} = q-1$, 
and $CSD(\lambda (b^\omega)) = CSD(b^\omega) - 8\pi^2 k \ell$, 
see Lemma \ref{index:shift}.  
 
When non-empty, $\hat\M^0(a, b)$ is isomorphic 
to $\hat\M (a^\omega, \lambda (b^\omega))$. We can prove that 
$\hat\M (a^\omega, \lambda (b^\omega))$ is empty, since the CSD functional 
is non-increasing along the gradient flow lines and  
the difference of CSD between $a^\omega$ 
and $\lambda (b^\omega)$ is negative: 
\[ 
\begin{array}{lll} 
&& CSD (a^\omega) -CSD (\lambda (b^\omega)) \\[2mm] 
&=&  CSD (a^\omega) -CSD (b^\omega) - (CSD  (\lambda (b^\omega))  
- CSD (b^\omega) ) \\[2mm] 
&=& CSD (a^\omega) -CSD (b^\omega) + 8\pi^2 k \ell <0, 
\end{array} 
\] 
as $k<0$ and $|CSD (a^\omega) -CSD (b^\omega)| < 8\pi^2\ell$.  
This proves that the entries below the diagonal are always zero.  
 
For $a\in \M_{Y_0, q}^{(\omega)}(\s)$ and  
$b\in \M_{Y_0, q-1}^{(\omega)}(\s)$, it is easy to 
see that we have
\[ 
\la \partial_n (a), b \ra = \la \partial_q^{(\omega)} (a), b\ra. 
\] 
 
For $a\in \M_{Y_0, q}^{(\omega)}(\s)$ and  
$b\in \M_{Y_0, q-1 + k\ell}^{(\omega)}(\s)$ with $ (k>0)$, 
we can define  
\[ 
\la \partial_{q, k}^{(\omega)} (a), b\ra = \la \partial_n(a), b\ra  
=\#\bigl( \hat\M^0(a, b)\bigr). 
\] 
This counts the points in the zero dimensional components
$\hat\M^0(a, b)$ in the moduli space of trajectories on $\B$  from $a$
to $b$. Equivalently, we have 
\[ 
\la \partial_{q, k}^{(\omega)} (a), b\ra = 
\#\bigl( \hat\M (a^{\omega}, \lambda(b^{\omega})) \bigr), 
\] 
where $\hat\M (a^{\omega}, \lambda(b^{\omega}))$ 
is  the moduli space of trajectories on $\B_\ell$  from  
$a^{\omega}$ to $\lambda(b^{\omega})$ and  
$\lambda$ represents the unique element $[\lambda]$ in  
$\G/\G_\ell$ such that  
\[ 
\la [\lambda] \cup  c_1(L), [Y_0]\ra = -kl. 
\] 
Thus we have $\ind_{Y_0} (\lambda(b^{\omega})) = q-1$. Therefore, if
non-empty, the moduli space $\hat\M (a^{\omega}, \lambda(b^{\omega}))$
is an oriented, compact 0-dimensional manifold with energy 
given by  
\[ 
8\pi^2 k\ell + CSD(a^{\omega}) - CSD(b^{\omega}) > 8\pi^2(k-1)\ell >0. 
\] 
This completes the proof of the Lemma. 
\endproof

The expression of $\partial_n$  and the appearance of  
$\partial_{q, k}^{(\omega)}$ in Lemma \ref{boundary:operator} 
lead us naturally to introduce a filtration of  
of the $\Z_\ell$ graded complex 
$C_*(Y_0,\s)$. The filtration is given by  
\[ 
F_{q}^{(\omega)} C_n =\bigoplus_{k \geq 0}C_{q+k\ell}^{(\omega)}(Y_0, \s).  
\] 
for $n \in \Z_\ell$, and $q\in \Z$ with $q\equiv n(mod \ \ell)$. 
Thus, we have 
\ba 
\cdots \subset F_{q + \ell}^{(\omega)} C_n 
\subset F_{q}^{(\omega)} C_n \subset F_{q -\ell}^{(\omega)} C_n 
\subset \cdots \subset C_n (Y_0, \s), 
\label{filtration} 
\na 
a finite length decreasing filtration of the $\Z$-graded 
Seiberg--Witten--Floer chain complex. From Lemma \ref{boundary:operator}, 
we see that  the boundary operator  
\[ 
\partial_n: \qquad F_{q}^{(\omega)} C_n \longrightarrow 
 F_{q-1}^{(\omega)} C_{n-1} 
\] 
preserves the filtration. Let  
$F_q^{(\omega)}H_n$ denote the homology of the complex 
\[ 
\cdots \stackrel{\partial}{\rightarrow}  
 F_{q}^{(\omega)} C_n \stackrel{\partial}{\longrightarrow} 
 F_{q-1}^{(\omega)} C_{n-1} \stackrel{\partial}{\rightarrow}  
\cdots.\] 
 
Define  
\[ 
F_q^{(\omega)} HF_n^{SW}(Y_0, \s) = Im \bigl(F_q^{(\omega)}H_n \rightarrow  
HF_n^{SW}(Y_0, \s) \bigr). 
\] 
We thus obtain a bounded 
filtration on $ HF_n^{SW}(Y_0, \s)$, 
\[\begin{array}{c} 
\cdots \subset F_{q + \ell}^{(\omega)} HF_n^{SW}(Y_0, \s)  
\subset F_{q}^{(\omega)}  HF_n^{SW}(Y_0, \s)\\[2mm] 
\subset F_{q -\ell}^{(\omega)}  HF_n^{SW}(Y_0, \s) 
\subset \cdots \subset  HF_n^{SW}(Y_0, \s). 
\end{array} 
\] 
 
The standard procedure of constructing the spectral sequence for a
filtration \cite{Sp} gives the following theorem on the relation
between the $\Z$-graded and the $\Z_\ell$-graded homology 
groups. 
 
\begin{The} 
There exists a spectral sequence $(E^k_{q, n}(Y_0, \s), d^k)$ 
with  
\[ 
E^1_{q, n}(Y_0, \s) \cong HF^{SW}_{q, (\omega)}(Y_0, \s)\] 
for $n\in \Z_\ell$ and $q\in \Z$ with $q\equiv n (mod \ \ell)$.  
The higher differentials 
\[ 
d^k: \qquad E^k_{q, n}(Y_0, \s) \longrightarrow 
E^k_{q-1 +k \ell , n-1 }(Y_0, \s)  
\] 
are induced by the maps $\partial_{q, k}^{(\omega)}$ defined in  
Lemma \ref{boundary:operator}. Furthermore,  the 
spectral sequence $(E^k_{q, n}(Y_0, \s), d^k)$ 
converges to the $\Z_\ell$-graded homology groups 
$HF_n^{SW}(Y_0, \s)$.  
\label{spectral:sequence} 
\end{The} 

\begin{proof} 
By construction, the $\Z_\ell$-graded chain complex $C_n(Y_0, \s)$ 
has a bounded filtration (\ref{filtration}) with the associated 
graded complex given by 
  \[ 
F_q^{(\omega)}C_n/F_{q+\ell}^{(\omega)}C_n 
= C^{(\omega)}_q(Y_0, \s). 
\] 
Then by the standard technique of \cite{Sp}, we derive the existence
of a spectral sequence   
\[ 
(E^k_{q, n}(Y_0, \s), d^k) 
\] 
with $E^1_{q, n}(Y_0, \s) \cong HF^{SW}_{q, (\omega)}(Y_0, \s)$, and  
\[ 
\begin{array}{l} 
Z_{q, n}^k(Y_0, \s) = \{ a \in F_q^{(\omega)} C_n | \  \partial_n (a) \in 
     F_{q-1+k\ell}^{(\omega)} C_{n-1} \}\\[2mm] 
E_{q, n}^k(Y_0, \s) = Z_{q, n}^k(Y_0, \s)/ \bigl( Z_{q+\ell, n}^{k-1}(Y_0, \s) 
+ \partial_{n+1} Z^{k-1}_{q+1 -(k-1)\ell, n+1}(Y_0, \s)\bigr). 
\end{array} 
\] 
The higher differentials are induced by $\partial_n$. The expression 
of $\partial_n$, as discussed in Lemma \ref{boundary:operator}, tells 
us that the higher differentials  acting  on $E_{q, n}^k(Y_0, \s)$ 
are defined by $\partial_{q, k}^{(\omega)}$. 
\end{proof}

\section{$\Z$-graded homology for exact triangles}

In this section, we present a different construction which applies to
the special case where $Y_0$ is obtained
by Dehn surgery with framing $0$ 
on a knot $K$ in a $\Z$-homology 3-sphere $Y$. Again we assume that
the $\spinc$-structure is non-torsion. The torsion case will be
analyzed in the next Section.

Under our assumptions, $Y_0$ has the homology of $S^1\times S^2$.  
We denote by $Y_1$ the $\Z$-homology 3-sphere obtained by  
Dehn surgery on $K$ with framing 1. With minor  
modifications, the arguments in this subsection can be extended 
to the case of a knot $K$ representing  
a zero homology class in $H_1(Y, \Z)$ in a general 3-manifold $Y$. 
 
In \cite{CMW}, we have identified the (perturbed) irreducible Seiberg--Witten 
monopoles (as critical points of the perturbed CSD functional) 
on $(Y, g)$, with monopoles on $(Y_1, g_1)$ and on $Y_0$, with
\ba 
\M_{Y, g, \mu}^* \cong \M_{Y_1, g_1}^* \cup \bigcup_{\s} \M_{Y_0}(\s). 
\label{3dsw:monopole} 
\na 
The metric and perturbation on $Y$ have been carefully chosen to 
met the above identification.  In particular, the perturbation $\mu$
is constructed as in \cite{CMW} so as to ``simulate the effect of
surgery''. Namely, it is obtained as a deformation of the circle of
flat connections on the tubular neighborhood of the knot $\nu(K)$ in $Y$,
so that it approximates the union of the lines of flat connections on
$\nu(K)$ inside $Y_1$ and $Y_0$.

For each $\spinc$ structure $\s$ 
on $Y_0$, with non-empty moduli space $\M_{Y_0}(\s)$,  
we have an injective map 
\ba 
 j_\s: \qquad \M_{Y_0}(\s) \longrightarrow 
   \M_{Y, g, \mu} 
\label{j_s} 
\na 
defined by the identification (\ref{3dsw:monopole}).
 
\begin{Def} Assume that $c_1(\s)$ is a non-torsion class 
in $H^2(Y_0, \Z)$. We define the $\Z$-grading of elements in
$\M_{Y_0}(\s)$ as  
induced from the grading on $\M_{Y, g, \mu}$, 
\[ 
\ind_{Y_0}^{(Y)} (a) = \ind_Y (j_\s (a)) 
\] 
for any $a\in \M_{Y_0}(\s)$. 
\end{Def} 
 
Recall that the grading on $\M_{Y, g, \mu}$ is defined by the 
spectral flow of the Hessian operator of CSD along a path 
from the critical point in $\M_{Y, g, \mu}$ to the unique 
reducible point $\vartheta$ in $\M_{Y, g, \mu}$.  
The relative index $\ind_{Y_0} (a, b)$ of two points   
$a$ and $b$ in  $\M_{Y_0}(\s)$, defined as the spectral flow of the Hessian 
operator of CSD function for $(Y_0, \s)$, agrees with the 
relative index $\ind_Y(j_\s(a), j_\s(b))$, cf. \cite{CMW}.  
Thus, we have 
\[ 
\ind_{Y_0}^{(Y)} (a) -\ind_{Y_0}^{(Y)} (b) 
= \ind_{Y_0} (a, b) \ (mod \ \ell) 
\] 
where $\ell \in \Z$ satisfies 
$c_1(\s) (H_2(Y_0, \Z)) = \ell \Z$. Thus, $\ind_{Y_0}^{(Y)} $ 
is a $\Z$-lifting of the $\Z_\ell$-grading on 
$\M_{Y_0}(\s)$.

\begin{Def} 
For any $q\in \Z$, we define 
\[ 
 \M_{Y_0}^{(q)}(\s)= 
 \{ a\in \M_{Y_0}(\s) \ | \ \ind_{Y_0}^{(Y)}(a)=q \}.  
\] 
The $q$-chains in the corresponding 
 $\Z$-graded Floer complex are the elements of 
the abelian group  
$ C_{(q)}(Y_0,\s) $  
generated by the points of $\M_{Y_0}^{(q)}(\s)$. 
The boundary operator 
\[ 
\partial_{(q)} : C_{(q)}(Y_0,\s) \to C_{(q-1)}(Y_0,\s)  
\] 
is defined as 
\[ 
\partial_{(q)}(a) = \sum_{b\in \M_{Y_0}^{(q-1)}(\s)} 
\# \bigl( \hat \M^0 (a, b)\bigr) b, 
\] 
where $\hat \M^0 (a, b)$ is the zero dimensional components of the 
moduli space of   
unparameterized flow lines in $\hat\M(a, b)$. The compactness theorem 
tells us that $\hat \M^0 (a, b)$ is an oriented, compact 0-dimensional
manifold. Thus, the coefficient $\#(\hat \M^0 (a, b))$ is well-defined 
as the algebraic sum of points in $\hat \M^0 (a, b)$.  
\label{(Y_0:q)} 
\end{Def} 
 
In order to prove that we have $\partial_{(q)} \circ \partial_{(q-1)}
=0$, we need the following Lemma to know
 when the zero-dimensional moduli space of flow lines between $a$
and $b$ over 
$(Y_0(r), \s)$ is non-empty. Here we follow the notation of
\cite{CMW}, where $Y_0(r)$ and $Y(r)$ denote the 3-manifolds $Y_0$ and
$Y$ endowed with long cylinders $T^2\times [-r,r]$ along the boundary
of the tubular neighborhood 
$\nu(K)$ around the knot $K$.
 
\begin{Lem} Let $\s$ be a $\spinc$ structure on $Y_0$ with 
non-empty $\M_{Y_0}(\s)$, and $c_1(\s)$ is a non-trivial class 
in $H^2(Y_0, \Q)$. Then the moduli spaces of flow lines can be 
described as follows. 
 Assume that $a, b$ are two critical points in $\M_{Y_0}(\s)$ 
of relative index 1 modulo $\ell$,
here $\ell \in \Z$ satisfies 
$c_1(\s) (H_2(Y_0, \Z)) = \ell \Z$.
If  for sufficiently 
large $r$, the $1$-dimensional component $\M_{Y_0(r)\times \R}^1(a,
b)$ of $\M_{Y_0(r)\times \R}(a, b)$, is non-empty, then we have
\[ 
 \ind_{Y_0}^{(Y)}(a) - \ind_{Y_0}^{(Y)}(b) =1. 
\]   
  \label{moduli:=} 
\end{Lem} 

\begin{proof} 
The result is  the consequence of the results on the geometric limits  
 in \cite{MW2} \cite{MW3}, together with the  dimension formulae
for various moduli spaces.

Suppose given a solution $[A_r(t), \Psi_r(t)]$ in a 
one-dimensional component of  
$\M_{Y_0(r)\times \R}(a, b)$, with asymptotic values
\[ 
[A_r(\pm \infty), \Psi_r(\pm \infty) ] = 
[A'_{\pm \infty}, \Psi'_{\pm \infty} ] \#_{a_{\pm}}^r [a^\pm, 0], 
\] 
for $t\to \pm \infty$. The geometric limits as $r \to
\infty$ are determined in \cite{MW2}.
Among these geometric limits, there are two parameterized
oriented  paths  $a_V(\tau), a_0(\tau)$ on 
$\partial_\infty (\M^*_{V})$ and  
$\M_{\nu(K)\subset Y_0} = L_{Y_0}$, respectively, where
both paths connect $[a^\pm, 0]$, consistently with the
orientations. Here $\M_{\nu(K)\in Y_0}$ is the
circle of flat $U(1)$-connections on $\nu(K)\subset Y_0$ modulo gauge
equivalence.  
There is also a holomorphic map  
$a_{D}$ from the  unit half-disc  to $\chi (T^2, Y_0)$
with boundary along   the paths $a_V(\tau)$ and $a_0(\tau)$.

First  we show that,  if the moduli space $\M_{Y(r)\times
\R}(j_\s(a), j_\s(b))$ is non-empty, then we have
  \[ 
 \M_{Y(r)\times \R}(j_\s(a), j_\s(b))  \cong 
 \M_{Y_0(r)\times \R}^1( a ,  b ). 
\] 
Under this assumption, the geometric limits of a family of solutions $[A_r(t),
\Psi_r(t)]$ in $\M_{Y_0(r)\times \R}(a,b)$, for $r \to \infty$, can
be deformed by a homotopy in $H^1(T^2,i\R)$ to geometric limits for a
monopole in $\M_{Y(r)\times \R}(j_\s(a),j_\s(b))$. 
In particular, this can happen if and only if
the path $a_0(\tau)=[A(\tau), 0]$ along
$L_{Y_0}$ can be homotopically deformed to a path along the curve  
$\partial_\infty \M_{\nu(K),\mu}$, by a homotopy that moves the
endpoints $a^\pm$ to the corresponding points $a^\pm(\epsilon)$,
along $a_V(\tau)$.
(We
use here the same notation as in 
\cite{CMW}, \cite{MW2}, \cite{MW3}.) Both paths $a_V(\tau)$ and 
$a_0(\tau)$ must be contractible
in $\chi (T^2, Y_0) = S^1 \times \R$.
In fact, the geometric limits of any
solution in $\M_{Y(r)\times
\R}(j_\s(a),j_\s(b))$ define an approximate solution to the monopole
equations on $Y_0(r)\times \R$, which can be deformed to a solution in 
the minimal energy component of 
$\M_{Y_0(r)\times \R}(a,b)$. Conversely, in this case, 
 the geometric limits of
$[A_r(t),\Psi_r(t)]$ can be spliced together to form a
solution to the monopole equations on $(Y(r)\times \R,\mu,g)$ between
the critical points $j_\s(a)$ and $j_\s(b)$. Thus the identification
of the moduli spaces implies that
\[ 
 \ind_{Y_0}^{(Y)}(a) - \ind_{Y_0}^{(Y)}(b) =1. 
\] 

Now assume that  the moduli space $\M_{Y(r)\times
\R}(j_\s(a), j_\s(b))$ is empty. Notice
that oriented paths $a_V(\tau), a_0(\tau)$
which are non-contractible in  the 
cylinder $\chi_0 (T^2, \nu(K)) \cong S^1 \times \R$ 
correspond to flowlines in the components with 
 higher  energy.   Here 
$$  \chi (T^2, Y_0)= \chi_0 (T^2, \nu(K))$$ 
is a covering of the torus
$\chi(T^2)$ of flat $U(1)$-connections on $T^2$ modulo gauge
equivalence, obtained by taking the  quotient only by those gauge
transformations on $T^2$ that extend to $Y_0(r)$, cf. \cite{CMW}.
Since we are considering the flowlines in the minimal energy
components $\M_{Y_0(r) \times \R}^{1}(a, b)$,
  the paths $a_V(\tau)$ and 
$a_0(\tau)$ are also contractible
in $\chi (T^2, Y_0) = S^1 \times \R$. The analysis
of the moduli space on the cobordant manifold $W_0$ 
with cylindrical ends $Y_0 (r) \times (-\infty, 0]$
and $Y(r) \times [0, \infty)$ in \cite{MW3} shows
that the holomorphic half-disc  $a_{D}$, 
from  the geometric limits of $[A_r(t), \Psi_r(t)]$ in the 
one-dimensional components of  
$\M_{Y_0(r)\times \R}(a, b)$, is the degenerate
limit of the assembled holomorphic triangles 
$\Delta^\epsilon$ as $\epsilon \to 0$, where $\epsilon$
is the parameter in the perturbation on $Y (r)$ that ``simulates
the effect of surgery". 
In particular, we have the following identification
\[
\M_{Y_0(r)\times \R}^1( a ,  b ) \cong
\bigcup_{c\in \M^*_{Y_1}} 
\M_{Y (r)\times \R} ( j_\s(a) , i(c) ) \times
\M ^{W_0, (0)} (i(c), b ).\]
Here $i(c)$ is the corresponding monopole
in $\M^*_{Y, \mu, g}$, under (\ref{3dsw:monopole}),
satisfying $\ind_Y (j_\s(a) ) - \ind_Y ( i(c) ) =1$,  and 
$\M ^{W_0, (0)} (i(c), b )$ is given by the non-empty 0-dimensional 
components of $\M ^{W_0, (0)} (i(c), b )$. Here we follow the notation
of \cite{MW3}. Then the dimension formulae developed in \cite{MW3} 
can be applied to show that 
\[
\ind_Y ( i(c) ) - \ind_{Y_0}^{(Y)}(b) =0.\]
Thus  we have 
$\ind_{Y_0}^{(Y)}(  a ) - \ind_{Y_0}^{Y}( b )=1$.   
This completes the proof of the statement. 
\end{proof}

Using the result of Lemma \ref{moduli:=}, we now show that
$\partial_{(q-1)}\circ \partial_{(q)} =0$, and we relate the resulting 
homology groups 
to the original $\Z_\ell$-graded Floer homology groups 
for $(Y_0, \s)$. 
 
\begin{Pro} 
\begin{enumerate} 
\item  For $n\in \Z_\ell$ and $q\in \Z$ with $q=n (mod \ \ell)$, then  
\[ 
C_n (Y_0, \s)=  \bigoplus_{k\in \Z}  C_{(q+ k \ell)}(Y_0,\s). 
\] 
\item Under the decomposition as 
above, the boundary operator $\partial_n$ on 
$C_n(Y_0, \s)$ can be expressed as  
$\oplus_{k\in \Z} \partial_{(q+k\ell)}$. Thus,
in particular, we have $\partial_{(q-1)}\circ 
\partial_{(q)} =0$ and  
\[ 
HF^{SW}_n(Y_0, \s) = \bigoplus_{k\in \Z} HF_{(q+k\ell)}^{SW}(Y_0, \s), 
\] 
where $HF_{(q)}^{SW}(Y_0, \s)$ is the q-th homology group 
of the chain complex $(C_{(q)}(Y_0, \s), \partial_{(q)})$. 
\end{enumerate}\end{Pro} 

\proof We only need to prove the second statement. 
For any generator $a\in C_{(q)}(Y_0, \s)$, as an element 
in $C_n (Y_0, \s)$, we shall verify that $\partial_n(a) = \partial_{(q)}(a)$. 
   From the definition, we know that 
\[ 
\partial_n(a) = \sum_{b\in \M_{Y_0}(\s): \ind_{Y_0}(a, b)=1(mod \ \ell)} 
\# \bigl( \hat \M^0_{Y_0\times \R} (a, b)\bigr) b,  
\] 
where $\hat \M^0_{Y_0\times \R} (a, b)$  
is the zero dimensional component of the 
moduli space of the  
unparameterized flow lines in $\hat\M _{Y_0\times \R}(a, b)$. 
By the result of Lemma \ref{moduli:=}, we know that the critical points 
$b \in  \M_{Y_0}(\s)$ with non-empty 
zero dimensional component in $ \hat \M_{Y_0\times \R} (a, b)$ 
are contained in the set  
\[ 
\M_{Y_0}^{(q-1)}(\s)= 
 \{ b\in \M_{Y_0}(\s) \ | \ \ind_{Y_0}^{(Y)}(b)=q-1 \}.  
\] 
This implies that $\partial_n (a) = \partial_{(q)} (a)$ for 
$a\in C_{(q)}(Y_0, \s)$, hence 
 $\partial_n  = \oplus_{k\in \Z} \partial_{(q+k\ell)}$.  
Then the statement follows from standard arguments  
in homological algebra. 
\endproof 
 
In the exact triangle for the Seiberg--Witten--Floer homologies  
relating $HF^{SW}_*(Y, g)$,  
$HF^{SW}_*(Y_1, g_1)$, and $HF^{SW}_{(*)}(Y_0, \s_k)$, as studied 
in \cite{MW2}\cite{MW3},   
the $\Z$-graded Seiberg--Witten--Floer homology for $(Y_0. \s)$ 
with $\s$ a non-trivial $\spinc$ structure, is 
precisely the  
$\Z$-graded $HF^{SW}_{(*)}(Y_0, \s)$ described in this subsection.

\begin{Rem}
In the case where $Y_0$ is obtained as 0-surgery on a knot $K$ in 
a homology 3-sphere, both the constructions of $\Z$-graded Floer
homology are possible. 
\end{Rem}

\section{Monopole homology for trivial $\spinc$-structures} 
 
In this section we consider again the general case where $Y_0$ is a
compact 3-manifold with $b_1(Y_0)>0$. We analyze the Seiberg--Witten
Floer theory in the case of a $\spinc$-structure $\s_0$ with 
$c_1(\det(\s_0))$ torsion. The results we obtain specialize to the
case of $Y_0$ a 0-surgery on a knot in a homology 3-sphere, with the 
trivial $\spinc$-structure with $c_1(\det(\s_0))= 0$.  

A cohomologically trivial perturbation 
would introduce a torus of reducible monopoles among the critical points 
of the Chern--Simons--Dirac functional given by the flat 
$U(1)$-connections $H^1(Y_0,\R)/H^1(Y_0,\Z)$. 
In order to avoid the torus of reducibles, we need to introduce a 
small cohomologically non-trivial perturbation of the Seiberg--Witten 
equations on $Y_0$ given by a 1-form $\rho$ satisfying 
$[*\rho] \in H^2(Y_0,\R)$ non-trivial. 
 With this perturbation,  $\M_{Y_0}(\s_0)$ 
consists of a finite set of points in $\B$, cut out transversely by 
the 3-dimensional (perturbed) Seiberg--Witten monopole equations.  
 
The different nature of the case $c_1(\s_0)$ torsion  
can be regarded as a result of the following phenomenon. In the case 
analyzed previously, with $c_1(\s)$ non-torsion, we studied the 
critical points of CSD on the infinite cyclic covering space $\B_\ell$ of 
$\B$. The main point of our argument has been the following: the non-trivial 
covering $\B_\ell$ is the essential tool in separating the various 
components of $\hat {\cal M}(a,b)$ of different dimensions,  
with each of them separately admitting a nice 
compactification to a smooth manifold with corners as in 
\cite{MW}, \cite{Mel}. 
 
In the case with $c_1(\s_0)$ torsion, however, the relative index 
$\ind (a)-\ind (b)$ is already well defined as an integer in
$\B$. Thus, the moduli space $\hat {\cal  
M}(a,b)$ does not break into components of different dimensions. In
fact, all components in $\hat \M (a, b)$ have constant dimension given 
by the relative index $\ind (a)-\ind (b) \in \Z$. 
In order to define the Floer homology, we still need a nice 
compactification of the moduli spaces $\hat {\cal M}(a,b)$. However, we 
lack the uniform energy bound on $\hat {\cal M}(a,b)$, since the CSD
functional is only circle-valued, due to the presence of the
non-trivial cohomology class of the perturbation $[*\rho]$.   
This lack of a uniform energy estimate implies that we cannot define
the usual boundary operator  
$$\la \partial a, b \ra= n_{ab} $$ 
of the Floer complex by setting  
$$ n_{ab}=\# \hat {\cal M}(a,b), $$ 
since $\hat{\cal M}(a,b)$ may be non-compact in this case.

Thus, we need to separate $\hat {\cal M}(a,b)$ into
components of   
uniform energy and define a boundary operator in terms of  
components of a fixed energy. We identify $\hat {\cal M}(a,b)$ 
with a subset of ${\cal M}(a,b)$ consisting of those trajectories 
with equal energy distributions on $(-\infty, 0]$ and $[0, \infty)$, that 
is, 
\ba\begin{array}{lll} 
&&\displaystyle{\int_{-\infty}^0} (\|\partial_t A(t)\|^2_{L^2(Y_0)} 
 + \|\partial_t \psi(t)\|^2_{L^2(Y_0)} ) dt \\[3mm] 
&=&  
\displaystyle{\int_0^\infty } (\|\partial_t A(t)\|^2_{L^2(Y_0)} 
 + \|\partial_t \psi(t)\|^2_{L^2(Y_0)} ) dt. \end{array} 
\label{energy:equal} 
\na 
 
Every non-constant 
 trajectory $[A(t), \psi(t)]$ in a non-empty connected 
 component of $\hat \M(a, b)$
has positive energy  
\[ 
\displaystyle{\int_\R} (\|\partial_t A(t)\|^2_{L^2(Y)} 
 + \|\partial_t \psi(t)\|^2_{L^2(Y)} ) dt > 0. 
\] 
The energy is constant, independent of the trajectory, in 
each connected component of $\hat \M(a, b)$. This energy agrees with  
the variation of the perturbed CSD functional along $[A(t),
\psi(t)]$. In fact, we have
$$ - \| \nabla CSD_{U,V} \|^2_{L^2} = \langle \partial_t (A,\psi), 
\nabla CSD_{U,V} \rangle, $$
if $CSD_{U,V}=CSD +U+V$ denotes the perturbed CSD functional, and 
$(A,\psi)$ is a solution of the flow equation
$$ \partial_t (A,\psi) = - \nabla CSD_{U,V}. $$
This gives
$$ \int_t^\infty \| \nabla CSD_{U,V}(A(\tau),\psi(\tau)) 
\|^2_{L^2(Y)} d\tau = $$
$$ CSD_{U,V}(A(t),\psi(t)) - CSD_{U,V}(A_b,\psi_b). $$

Thus, we can define the energy function over 
$\hat \M(a, b)$ as 
\ba 
\begin{array}{c} 
\E: \qquad \hat \M (a, b) \longrightarrow \R \\[2mm] 
\E([A(t), \psi(t)]) =  
\displaystyle{\int_\R} (\|\partial_t A(t)\|^2_{L^2(Y)} 
 + \|\partial_t \psi(t)\|^2_{L^2(Y)} ) dt. 
\end{array} 
\label{energy:E} 
\na 
 
The image of $\E$ is positive and discrete  in $\R$. The energies of
any two connected components may differ by 
some multiple of a positive number  
\[ 
e_\rho = min  \{ | [*\rho] ([\Sigma])|  : [\Sigma] \in H_2(Y_0, \Z) \}. 
\] 
 
Denote $e_{min}(a, b)$ the minimal energy on $\hat\M(a, b)$.  
Then, for any $n\ge 0$ in $\Z$,  
we can define  
\[ 
\hat\M^{(n)} (a, b) = \E^{-1} (e_{min}(a, b) + n e_\rho ). 
\] 
We have the following compactness theorem for the 
moduli space $\hat\M^{(n)}(a, b)$ with the fixed  
energy $e_{min}(a, b) + n e_\rho $, for some $n\ge 0$. 
 
\begin{Pro} 
Let $(Y_0, \s_0)$ be a 3-manifold with $b_1(Y)>0$ and with 
$c_1(\s_0)$ torsion. Let $\hat\M^{(n)}(a, b)$ be the component 
of $\hat \M(a, b)$ defined as above, for two critical points 
$a$ and $b$ of relative index $q+1$. Then  
$\hat\M^{(n)}(a, b)$ can be compactified to be a smooth manifold 
of dimension $q$ with corner structure, 
by adding lower dimensional boundary strata. 
The codimension $k$ boundary strata are given by  
\[ 
\bigcup_{a_1, \cdots, a_k} \hat\M^{(n_0)}(a, a_1) 
\times \hat\M^{(n_1)}(a_1, a_2) \times \cdots 
\times \hat\M^{(n_k)} (a_k, b). 
\] 
Here the union is over all possible sequences $a_0=a, a_1, \cdots 
a_k, a_{k+1} =b$ with strictly decreasing indices and with 
$\sum_{j=0}^k  n_j = n$.  
In particular, we have the following results. 
\begin{enumerate} 
\item The component of fixed energy in 
$\hat\M(a, b)$ with relative index $1$ is compact. 
\item If $a$ and $c$ are two critical points of relative index  
$\ind(a) -\ind(c) =2$, then
the compactified  moduli space  
$\hat\M^{(n)}(a,c)$ is an oriented, smooth manifolds with boundary 
points given by  
\[ 
\bigcup_{b, i+j =n} \hat\M^{(i)}(a, b) \times \hat\M^{(j)}(b, c), 
\] 
where $b$ runs over all critical points with relative 
index $\ind(a) -\ind(b) =1$.  
\end{enumerate} 
\label{compact:c_1=0} 
\end{Pro} 

\proof 
The convergence and gluing theorems developed  
in Section 4 of \cite{MW} can be applied to this case to give
$\hat\M^{(n)}(a, b)$ the structure of 
a smooth manifold structure with corners. In particular, the local
diffeomorphism provided by the convergence and gluing theorem in \cite{MW} 
is energy preserving. Then the above claims follow from 
the arguments of Section 4 of \cite{MW}. Notice that we have both
inequalities 
$$ e_{min}(a, b)\leq e_{min}(a,a_1) + \cdots + e_{min}(a_k,b) $$
and
$$ e_{min}(a, b)\geq e_{min}(a,a_1) + \cdots + e_{min}(a_k,b), $$
where we use essentially the identification between energy and
variation of the perturbed CSD functional.  
\endproof

We  now introduce  a chain  complex over the ring $\Z[[t]]$ of 
formal power series as follows. We define
the chain  complex $C_*(Y_0, \s_0, \Z[[t]])$ 
to be the $\Z[[t]]$-module generated by the critical points in
$\M_{Y_0}(\s_0)$,  
with the grading given by the relative grading $\ind(a)-\ind(a_0)$, where 
$a_0$ is a fixed critical point in $\M_{Y_0}(\s_0)$. The boundary 
operator $\partial$ is defined by 
\ba  
\partial_{[t]} ( a )= \sum_{b :\ind (a)-\ind (b)=1} \sum_{k\geq 0}  
 \#\bigl( \hat\M^{(n)} (a, b) \bigr) t^n b. 
\label{partial:[[t]]} 
\na  
 
By the result of Proposition \ref{compact:c_1=0}, we know that 
the boundary operator $\partial_{[t]}$ in  
 (\ref{partial:[[t]]}) is well-defined and 
 satisfies $\partial ^2=0$. 
The corresponding  homology groups are  the homology groups 
of the chain  complex $(C_*(Y_0,\s_0,\Z[[t]]), \partial_{[t]})$, 
\ba 
 HF_*(Y_0,\s_0,\Z[[t]])= H_*(C_*(Y_0,\s_0,\Z[[t]]), \partial_{[t]}).  
\label{SWF:[[t]]} 
\na   
 
In the construction of the exact triangle in \cite{CMW}, \cite{MW2},
\cite{MW3}, we set  
$$ HF_{(*)}(Y_0,\s_0)=HF_*(Y_0,\s_0,\Z[[t]])|_{t=0}.$$ 
Here the grading on $\M_{Y_0}(\s_0)$ is induced from the 
injective map 
\[ 
j_{\s_0}: \qquad  
 \M_{Y_0}(\s_0) \longrightarrow  \M_{Y, g, \mu}, 
\] 
that is, for any $a\in \M_{Y_0}(\s_0)$, we define 
the induced grading from $j_{\s_0}$ as 
\[ 
\ind_{Y_0}^{(Y)} (a) = \ind_{Y} (j_{\s_0}(a)). 
\] 
This induced grading on $\M_{Y_0}(\s_0)$ satisfies (cf. \cite{CMW}) 
\[ 
\ind_{Y_0}^{(Y)} (a) -\ind_{Y_0}^{(Y)} (b) 
= \ind_{Y_0} (a) - \ind_{Y_0} (b). 
\] 
After a suitable choice of a base point $a_0$ in $\M_{Y_0}(\s_0)$, 
we can assume that  
\ba 
\ind_{Y_0}^{(Y)} (a) = \ind_{Y_0}(a) = \ind_Y(j_{\s_0}(a)). 
\label{index:=} 
\na 

Notice that, in the case of the 0-surgery $Y_0$ that has
$b_1(Y_0)=1$, the Floer homology $HF_*(Y_0,\s_0,\Z[[t]])$ defined as in
(\ref{SWF:[[t]]}) depends on the cohomology class $[*\rho]\in
H^2(Y_0,\R)$. This dependence can be seen from the fact that, when the
value $[*\rho]\in \R$ crosses zero, some irreducible solutions in
$\M_{Y_0,\rho}(\s_0)$ may hit the circle of reducibles, cf. \cite{Lim}. 
 
We have the analog of Lemma \ref{moduli:=}, in the case of a
3-manifold $Y_0$ that is obtained as a 0-surgery on a knot in a
homology 3-sphere $Y$. 
 
\begin{Lem} Let $\s_0$ be a trivial $\spinc$ structure on $Y_0$, where
$Y_0$ is obtained as a 0-surgery on a knot in a
homology 3-sphere $Y$.  Let $a, b$ be two critical points 
in $\M_{Y_0}(\s_0)$ with relative index $ 1$. Then 
if  $ \hat\M_{Y_0(r)\times \R}^{(0)} (a, b)$ is
non-empty, for sufficiently large $r$
the oriented paths $a_0(\tau)$, $a_V(\tau)$ in the geometric 
limits of any flow line in the component of minimal energy 
are contractible in $\chi_0(T^2,Y_0)$.
Here the moduli spaces $\M_{Y_0}(\s_0)$
and $\hat\M_{Y_0(r)\times \R}^{(0)} (a, b)$ are considered with
an additional perturbation $\rho$ with $[*\rho]=\eta >0 $ in
$H^2(Y_0,\R)=\R$,  as in \cite{CMW}, cf. Remark
\ref{epsilon:rem}. 
\label{s_0moduli:=} 
\end{Lem} 
\proof 
This is again a direct applications of the geometric limit results 
 in \cite{MW2} and \cite{MW3}. The argument proceeds as in the proof 
of  Lemma \ref{moduli:=}. 
 
\endproof 
 
The homology groups 
for $(Y_0, \s_0)$,  which appear in the exact triangle, 
are precisely given by the homology groups 
$$ HF_{(*)}(Y_0,\s_0)=HF_*(Y_0,\s_0,\Z[[t]])|_{t=0},$$ 
after a possible global index shifting.

\begin{Rem} The  Seiberg--Witten--Floer homology for $(Y_0, \s_0)$ 
depends on the perturbation $[*\rho]$, as the energy function  
$\E$ and the Novikov-type chain complex all depend on the 
choice of $[*\rho]$. In the exact triangle, we use  
a particular choice of $[*\rho]$, satisfying $$ [*\rho]
=\eta PD_{Y_0}(m) \in H^2(Y_0,\R), $$  
with $m$ the meridian of the knot $K$.
\label{epsilon:rem}
\end{Rem}

\section{Gluing formulae along  the trivial $\spinc$ structure} 
 
Motivated by the general ideas of topological field 
theory, one can associate to any 4-manifold $(X_i, \s_i)$ with a
boundary $(Y_0, \t)$ Seiberg--Witten invariants with values in
the monopole homology groups $HF^{SW}_*(Y_0, \t)$,
\[ 
SW_{X_i}(\s_i, \cdot ): \qquad \AAA (X_i) \longrightarrow HF^{SW}_*(Y_0, \t), 
\] 
where $\AAA (X_i) = Sym^*(H_0(X_i)) \otimes \Lambda^*(H_1(X_i))$ 
is the free graded algebra generated by 
elements in $H_0(X_i)$ of degree $2$ and cycles in $H_1(X_i)$ 
of degree $1$.  
Moreover,  
if a closed 4-manifold $X=X_1\cup_{Y_0} X_2$ splits in two components
$X_i$, $i=1,2$, along a 3-manifold $Y_0$,  
then the natural 
pairing between $HF^{SW}_*(Y_0, \t)$ and $HF^{SW}_*(-Y_0, -\t)$ 
can be applied to give gluing formulae for the 
Seiberg--Witten invariants of $X$ in terms of the $SW_{X_i}$. 
For any 3-manifold $(Y_0, \t)$ with $b_1(Y)>0$ and $c_1(\t)$
non-torsion, these gluing formulae were developed in  
\cite{CW}.  For a manifold $Y$ with a torsion $\spinc$ structure $\t$, 
we discuss the corresponding gluing formulae in this section. 
 
To the purpose of gluing relative invariants, it is natural to  
construct the monopole homology groups over the  
field $\Q((t))$ of formal Laurent series and use  
another definition of the energy in the construction of 
monopole homology groups. 
 
We assume that the perturbation term 
 $[*\rho]$ is Poincar\'e dual to $\epsilon [\Gamma]$ for 
some 1-cycle $\Gamma$ in $Y_0$ and a small $\epsilon \in \R^+$. Here
$\Gamma$ is an integer cycle, that is, it defines a (non-torsion) integer
homology class.  The moduli space of critical points 
with this perturbation is denoted by  
$\M_{Y_0}(\t, \Gamma)$. 
Let $a, b$ be two critical points in $\M_{Y_0}(\t, \Gamma)$ 
 of relative index $k+1>0$. We know 
that $\hat\M (a, b)$ has infinite many components of dimension $k$.  
 Note that, if we identify elements in  
$\hat\M (a, b)$ with trajectories in $\M (a, b)$ satisfying 
the equal energy condition (\ref{energy:equal}), we can now  
define the following map on $\hat\M (a, b)$: 
\ba 
\begin{array}{c} 
\E_\rho: \qquad \hat\M (a, b) \longrightarrow \Z\\[2mm] 
\E_\rho ([A(t), \psi(t)]) = \la
[\displaystyle{\frac{\sqrt{-1}}{2\pi}}F_{A(t)}],  
[\Gamma \times \R]\ra. 
\end{array} 
\na 
Here $F_{A(t)}$ is the 4-dimensional curvature of $A(t)$. 
$\E_\rho$ essentially measures the variation of the Chern-Simons-Dirac 
functional on $Y_0\times \R$ up to a positive scalar.  
We denote by
\[ 
\hat\M^{[n]}(a, b) = \E^{-1}_\rho(n) 
\] 
the preimage of $n$ under the energy map $\E_\rho$. 
 
\begin{Lem} Let $Y$ be a 3-manifold with $b_1(Y)>0$, endowed with a
$\spinc$ structure $\t$ with $c_1(\t)=0$.
Let $a, b$ be two critical points in $\M_{Y_0}(\t, \Gamma)$ 
 of relative index $k+1>0$. We have the following results.
\begin{enumerate} 
\item $\hat\M^{[n]}(a, b)$ is empty for $n<<0$; 
\item For any non-empty $\hat\M^{[n]}(a, b)$, there exists 
a compactification of $\hat\M^{[n]}(a, b)$ to a smooth manifold with
corners, obtained by adding lower dimensional strata of the   
form  
\[ 
\bigcup_{a_1, \cdots, a_k} \hat\M^{[n_0]}(a, a_1) 
\times \hat\M^{[n_1]}(a_1, a_2) \times \cdots 
\times \hat\M^{[n_k]} (a_k, b). 
\] 
Here the union is over all possible sequences $a_0=a, a_1, \cdots 
a_k, a_{k+1} =b$ in $\M_{Y_0}(\t, \Gamma)$ 
 with strictly decreasing indices and with 
$\sum_{j=0}^k  n_j = n$.  
\end{enumerate} 
In particular, given any two critical points $a$ and $b$ of relative
index $1$, if the moduli space $\hat\M^{[n]}(a, b)$ is non-empty, then 
it is compact.  For  
any two  critical points $a$ and $c$ of relative index $2$,  
the compactified  moduli space  
$\hat\M^{[n]}(a,c)$ is an oriented, smooth manifold with boundary 
points given by  
\[ 
\bigcup_{b, i+j =n} \hat\M^{[i]}(a, b) \times \hat\M^{[j]}(b, c). 
\] 
Here $b$ runs over all critical points with relative 
index $\ind(a) -\ind(b) =1$.  
\label{compact:[n]} 
\end{Lem} 

\proof 
Any non-empty $\hat\M^{[n]}(a, b)$ has positive and constant  
energy defined by $\E$  (\ref{energy:E}). This implies 
that $\hat\M^{[n]}(a, b)$ must be empty for $n<<0$. The argument for
the compactification is analogous to
Proposition \ref{compact:c_1=0}. 
\endproof 
 
\begin{Def} Define the chain complex  $C_{\Gamma, *}(Y_0, \t,
\Q((t)))$ as the vector space generated by the critical points in
$\M_{Y_0}(\t, \Gamma)$, over  
the coefficient field of formal Laurent series $\Q((t))=\Q[[t]][t^{-1}]$, 
with the grading given by the relative index. The boundary operator  
is given by 
\[ 
\partial (a) = \sum_{b: \ind_{Y_0}(a) -\ind_{Y_0} (b) =1} 
   \sum_{n\in Z} \#\bigl( \hat\M^{[n]} (a, b) \bigr) t^n b. 
\] 
The fact that $\partial^2 =0$ follows from Lemma \ref{compact:[n]}. 
The homology groups of the chain complex $C_*(Y_0, \t, \Q((t)))$ 
are denoted by $HF^{SW}_{\Gamma, *}(Y_0, \t, \Q((t)))$.  
\end{Def} 
 
By the construction of $HF^{SW}_{\Gamma, *}(Y_0, \t, \Q((t)))$, we
know that these homology groups 
depend on the choice of the non-trivial cohomology 
class $[*\rho]$, which in turn depends on the choice of the 
one-cycle $\Gamma$.  There exists a natural isomorphism  
$HF^{SW}_{\Gamma, *}(Y_0, \t, \Q((t)))  
\cong HF^{SW, *}_{-\Gamma}(-Y_0, -\t, \Q((t)))$, where 
$-Y_0$ is the manifold $Y_0$ with the reversed orientation. For the
definition of Seiberg--Witten--Floer cohomology, in the 
equivariant setting, the reader can refer to Section 5 of \cite{MW}.
This yields a natural pairing 
\ba 
\la, \ra:  HF^{SW}_{\Gamma, *}(Y_0, \t, \Q((t)))  
\times HF^{SW}_{-\Gamma, *}(-Y_0, -\t, \Q((t))) 
\longrightarrow \Q((t)). 
\label{pairing} 
\na 
 
\subsection{Relative Seiberg--Witten invariants and gluing formulae} 
 
Let $X_1$ be an oriented Riemannian  4-manifold with a cylindrical 
end modeled on  $Y_0 \times [0, \infty )$. Let $\s_1$ be a $\spinc$
structure on $X_1$, such that the induced $\spinc$ structure $\t$ on
$Y_0$ satisfies 
$c_1(\t)=0 $ in $H_2(Y_0, \Q)$. Denote by $\M_{Y_0} (\t, \Gamma)$ 
 the moduli space of   
Seiberg--Witten monopoles on $(Y_0, \t)$, with respect to a perturbation  
$\rho$ satisfying $[*\rho]= \epsilon PD([\Gamma])$, for a non-trivial 
1-cycle $\Gamma$ in $Y_0$. Here again we assume that $\Gamma$ defines
a non-torsion integer homology class.
Moreover, for simplicity, we assume that $\Gamma$ lies
in the image of the map $H_2(X_1, Y_0; \Z) \rightarrow H_1(Y, \Z)$.
As we discuss in this section, the choice of $\Sigma$ with 
$\partial \Sigma = \Gamma$ will provide  a convenient way of 
separating the moduli spaces on $(X_1, \s_1)$ into
components of uniform energy, which admit a nice compactification
by adding lower dimensional strata.

 We assume also that we are 
working with an 
additional generic perturbation $(U, V)$ as in section 2, in order to
achieve transversality for the moduli space 
of trajectories.   We write 
$\eta_0=  \rho + \sum_j \frac{\partial U}{\partial \tau_j}\mu_j$ and 
$\nu_0 = \sum_j \frac{\partial V}{\partial \zeta_j}\nu_j$. 
Choose $c: X_1 \to \R^+$ to be 
a cut-off function supported  away from a compact set and equaling $1$ 
on the end. Choose $\mu \in  
 \Lambda ^{2, +}(X_1, i\R) $, satisfying $\mu - c\sqrt{-1}\eta_0^+ \in
C_\delta^k$. Here $C_\delta^k$ denotes the space of
$\delta$-decaying $C^k$-forms for some large integer $k$ and $\eta_0^+ $ 
is the self-dual 2-form obtained by $\eta_0$.  
Next we choose an imaginary valued 1-form $\nu$ on $X_1$, with 
$\nu - c\nu_0  \in C_\delta^k$.  Over the cylindrical end, we write 
$(A, \psi)$ in temporal gauge. When restricted to  
the cylindrical end, $\mu$ and $\nu$ can be thought of
as functions of $(A, \psi)$ on the cylindrical end. 
The relative Seiberg--Witten invariant for $X_1$ is obtained from
the perturbed Seiberg--Witten moduli space on $X_1$ with asymptotic limit 
representing a critical point in $\M_{Y_0} (\t, \Gamma)$.   
 
For each $a$ in $\M_{Y_0} (\t, \Gamma)$, define 
$\M_{X_1}(\s, a)$  to be the moduli space of 
the Seiberg--Witten equations on $X_1$ with asymptotic  
limit in the class of $a$, modulo the group of gauge transformations. 
The function spaces for the monopole equations modeled on $a$ are 
$L_1^2$-Sobolev spaces with `$\delta$-decay' where $\delta > 0$ is 
determined by  the least absolute eigenvalue of the Hessian of 
the CSD at $a$. Let $(A_a, \psi_a)$ be a  
$\spinc$ connection and a spinor on $X_1$ which agrees with  
the pull-back of some gauge 
representatives of $a$  on the cylindrical end. Note that the 
gauge class of  
$(A_a, \psi_a)$ can be identified with the set 
$H^1(Y_0, \Z)/ Im (H^1(X_1, \Z)\to  H^1(Y_0, \Z))$. 
 
Let $\M_{X_1}(\s_1, a)$ be the moduli space of solutions 
to the perturbed Seiberg--Witten equations on $(X_1, \s)$: 
\ba 
\left\{\begin{array}{l} 
F_A^+= q (\psi )+ \mu  \\[2mm] 
D_A \psi + \nu. \psi = 0, 
\end{array}\right. 
\label{SW:cyl} 
\na 
for $(A, \psi)$ in the weighted Sobolev space 
\[ 
\{ (A, \psi)| A- A_a \in L_{k, \delta }^2, 
                         \psi -\psi_a \in  L_{l, \delta, A_a}^2\}, 
\] 
and with gauge group given by $\{ u : X \to \C | 
 |u| =1 , 1-u \in L_{l+1, \delta}^2\}$. 
Here we fix a $U(1)$-connection $A_a$ on the spinor bundle to define the 
covariant derivatives on the spinor sections. 
For generic $(\mu, \nu)$, the moduli space $\M_{X_1}(\s_1, a)$, 
if non-empty, is a smooth, oriented manifold 
whose dimension is given by the index of the deformation complex 
for (\ref{SW:cyl}). Denote this index by $\ind_{X_1}(a)$. We then 
have  
\[ 
\ind_{X_1}(a) - \ind_{X_1}(b) = \ind_{Y_0}(a) -\ind_{Y_0}(b). 
\] 
In general, we cannot compactify $\M_{X_1}(\s_1, a)$ 
naturally by adding lower dimensional strata.  
Instead, we define some energy map on $\M_{X_1}(\s_1, a)$. 
Assume that the 1-cycle $\Gamma$ is bounded by a 2-cycle $\Sigma$ 
in $X_1$ with $\partial \Sigma = \Gamma$. Then the energy 
map is defined as the relative Chern class of  
$[A, \psi]$ restricted to $(\Sigma, \Gamma)$: 
\ba 
\begin{array}{c} 
\E_\Sigma: \qquad \M_{X_1}(\s_1, a) \longrightarrow 
 H^2(\Sigma, \Gamma; \Z)\\[2mm] 
\E_\Sigma ([A, \psi]) = \la [\displaystyle{\frac {\sqrt{-1}}{2\pi}} 
F_A] , [\Sigma]\ra. 
\end{array} 
\label{energy:Sigma} 
\na 
Denote by $\M_{X_1}^{[n]} (\s_1, a)$ the preimage of $n\in \Z$ under the 
energy map $\E_\Sigma$. Then we can compactify $\M_{X_1}^{[n]} (\s_1, a)$ 
to a smooth manifold with 
corners. The following proposition can be proved by adapting  the 
convergence and gluing theorems of \cite{MW} to the present setting. 
 
\begin{Pro} \label{compactify:X_1} 
Assume that  the moduli space $\M_{X_1}(\s_1, a)$ is non-empty. Then 
$\M_{X_1}^{[n]} (\s_1, a)$ is empty for $n<< 0$. Moreover, if non-empty,
$\M_{X_1}^{[n]} (\s_1, a)$ has a compactification to a smooth manifold 
with corners, with lower dimensional strata of the  
form  
\[ 
\bigcup_{a_1, \cdots, a_k} \M_{X_1}^{[n_0]} (\s_1, a_1)  
\times \hat\M^{[n_1]}(a_1, a_2) \times \cdots 
\times \hat\M^{[n_k]} (a_k, a). 
\] 
Here the union is over all possible sequences $ a_1, \cdots 
a_k, a_{k+1} =a$ in $\M_{Y_0}(\t, \Gamma)$ 
 with strictly decreasing indices and  with
$\sum_{j=0}^k  n_j = n$.  
\end{Pro} 

\proof  From the Weitzenb\"ock formula and the Seiberg--Witten equations,
we know that $\E_\Sigma$ is bounded from below. Thus, 
$\M_{X_1}^{[n]} (\s_1, a)$ is empty for $n<<0$. 
The remaining statements follow from the convergence and 
gluing theorems of \cite{MW}. 
\endproof 
 
Denote by $\AAA(X_1)$ the free graded algebra 
\[ 
Sym^*(H_0(X_1)) \otimes \Lambda^* (H_1(X_1))  
\] 
where the degree of the elements 
in $H_0(X_1)$ is 2 and the degree of the elements 
in $H_1(X_1)$ is 1. For any monomial $z = [x]^n \gamma_1 \wedge
\gamma_2 \wedge \cdots  
\wedge \gamma_k$ in $\AAA(X_1)$, and for each $m\in \Z$, 
 we need to consider all the  
components of dimension $2n +k$ in  
$\M_{X_1}^{[m]} ( \s_1) = \cup_a \M_{X_1}^{[m]}(\s_1, a)$. 
Suppose that $\M ^{[m]}_{X_1}(\s_1,  a_{(2n+k)})$ is the union of components 
 of dimension $2n+k$ and has energy $m$ 
under the energy map (\ref{energy:Sigma}). Here the
$\{a_{(2n+k)}\}$ denote all  those 3-dimensional monopoles 
(with perturbations) on $(Y_0, \t)$ which satisfy 
\[ 
\ind_{X_1} (a_{(2n+k)} ) = 2n+k. 
\] 
Thus, any non-empty component of $\M_{X_1}(\s_1,  a_{(2n+k)})$ 
has dimension $2n+k$.  
 Choose smooth loops 
which represent 1-cycles $\gamma_1, \cdots, \gamma_k$. 
The holonomy map of the  
$[A]$-part in $[A, \psi]$ around these   
$k$ loops defines maps 
\[ 
hol: \qquad \M  ^{[m]}_{X_1}(\s_1,  a_{(2n+k)}) \mapsto U(1) ^k. 
\] 
These maps can be extended naturally to the compactification of the
moduli spaces
$\M  ^{[m]}_{X_1}(\s_1,  a_{(2n+k)})$, as described by Proposition 
\ref{compactify:X_1}. The  generic fiber of $hol$ defines a smooth 
submanifold $V_{\gamma_1, \cdots, \gamma_k}$ of dimension $2n$ 
in $\M ^{[m]}_{X_1}(\s_1,  a_{(2n+k)}) $. This submanifold 
  $V_{\gamma_1, \cdots, \gamma_k}$ can be 
compactified by the corresponding fiber of $hol$ on the 
boundary strata of $\M  ^{[m]}_{X_1}(\s_1,  a_{(2n+k)})$. 
 
On $\M ^{[m]}_{X_1}(\s_1,  a_{(2n+k)})$, there is 
a canonical $U(1)$  bundle, defined by the based monopole moduli space 
$\tilde \M ^{[m]}_{X_1}(\s_1,  a_{(2n+k)})$,  
which consists of those solutions representing points in  
$\M ^{[m]}_{X_1}(\s_1,  a_{(2n+k)})$ but considered modulo only
those gauge transformations which fix the fiber of the spinor bundle  
over $x$. A different choice of base point provides an isomorphic 
$U(1)$ bundle.  Therefore, we have an associated complex 
rank $n$ vector bundle over $\M  ^{[m]}_{X_1}(\s_1,  a_{(2n+k)})$ 
\[ 
 E_n (  a_{(2n+k)}) = \tilde \M ^{[m]}_{X_1}(\s_1,  a_{(2n+k)}) \times
_{U(1)} \C^n.  
\] 
This vector bundle $  E_n( a_{(2n+k)})$ is compatible with the
boundary strata   
of $ \M ^{[m]}_{X_1}(\s_1,  a_{(2n+k)})$ in the sense that $ E_n (a_{(2n+k)})$ 
has a natural extension to the boundary strata of  
$ \M ^{[m]}_{X_1}(\s_1,  a_{(2n+k)})$, to a vector bundle in the
category of manifolds with corners (cf. \cite{Mel}). 
Note that $  E_n( a_{(2n+k)})$ is isomorphic to  
$\oplus_{i=1}^n {\cal L}(x_i)$, where ${\cal L}(x_i)$ is 
complex line bundle associated with the $U(1)$-bundle  
$\tilde \M ^{[m]}_{X_1}(\s_1,  a_{(2n+k)})$  
with the base point at $x_i \in X_1$. By construction, the line bundle
${\cal L}(x_i)$ extends to a line bundle in
the category of manifolds with corners over the
compactification. Since $\oplus$ is well defined in the category of
manifolds with corners, so does $E_n( a_{(2n+k)})$.
 
Now consider the restriction of $E_n (a_{(2n+k)})$ to the   
$2n$-dimensional submanifold $V_{\gamma_1, \cdots, \gamma_k}$, 
which is also compatible with its compactification. 
 
Choose a generic transversal section $\sigma_0$,  which 
is also transversal along the boundary strata, i.e. a strata
transverse section of a bundle over a manifold with corners.
By the transversality along the boundary strata, we know 
that all the zeroes of this section in $V_{\gamma_1, \cdots, \gamma_k}$ 
lie within a compact set and consist of finitely many points
with an orientation. 
Counting these points gives a number which is 
denoted by $SW_{X_1}^{[m]}(\s_1, z;  a_{(2n+k)})$. The relative
Seiberg--Witten invariant is defined to be 
\[ 
SW_{X_1}(\s_1, z) = \sum_{m\in \Z, a_{(deg (z))}} SW_{X_1}^{[m]} 
 (\s_1, z;  a_{(deg(z))})  < a_{(deg(z))}>t^m, 
\] 
as an element in the Seiberg--Witten--Floer complex. 
It is a routine to prove that this element is in fact a cycle, thus
defining an element in the Seiberg--Witten--Floer homology groups 
$HF_{\Gamma, *}^{SW}(Y_0, \t, \Q((t)))$. We also have the analog of
Proposition 3.3 in \cite{CW}, stating that  
$SW_{X_1}(\s_1, z)$, as an element in  
$HF_{\Gamma, *}^{SW}(Y_0, \t, \Q((t)))$, is independent of  the 
choice of $V_{\gamma_1, \cdots, \gamma_k}$ and the choice 
of the strata transversal section of $E_n$. Note that, if we choose
another relative 2-cycle $\Sigma$ in the definition of
the energy map (\ref{energy:Sigma}), the relative
invariant $SW_{X_1}(\s_1, z)$ only changes by  multiplication by a
certain power of $t$.

\begin{Pro} 
$SW_{X_1}(\s_1, z) \in HF_{\Gamma, *}^{SW}(Y_0, \t, \Q((t)))$ is well-defined, 
that is,
it is independent of the choices of $V_{\gamma_1, \cdots, \gamma_k}$
and the choice 
of the strata transversal section $\sigma_0$ of $E_n$.
\end{Pro}
\begin{proof}For any $m\in \Z$, let $V'_{\gamma_1, \cdots, \gamma_k}$
be another  
$2n$-dimensional submanifold in $\M^{[m]}_{X_1}(\s_1, a_{(2n+k)})$  
obtained
from the
construction of the 
holonomy map and let  
$\sigma_1$ be a strata
 transversal section of $E_n$ over 
$V'_{\gamma_1, \cdots, \gamma_k}$. We will prove
that the difference between  $\# \sigma_0^{-1} (0)$ and 
$\# \sigma_1^{-1}(0)$
defines an element in $C_{\Gamma, *}(Y_0, \t, \Q((t)))$
 which is homologous to zero.

Over $\M^{[m]}_{X_1}(\s_1, a_{(2n+k)})\times [0, 1]$, we can choose a
strata transverse section $\sigma$ of $ E_0$, 
with $\sigma(\cdot , 0) = \sigma_0$ and $\sigma(\cdot , 1) =\sigma_1$.
We also choose a 
codimension $k$ submanifold $V$ in 
$\M^{[m]}_{X_1}(\s_1, a_{(2n+k)}) \times [0, 1]$,  whose intersection
with $\M^{[m]}_{X_1}(\s_1,  a_{(2n+k)}) \times \{ 0\}$ is 
$V_{\gamma_1, \cdots, \gamma_k}$, and whose intersection with 
$\M^{[m]}_{X_1}(\s_1,  a_{(2n+k)}) \times \{ 1\}$ is 
$V'_{\gamma_1, \cdots, \gamma_k}$, respectively.  $V$ is also compatible
with the boundary strata of $\M^{[m]}_{X_1}(\s_1, a_{(2n+k)}) \times [0, 1]$.

   From this construction, the zero set of $\sigma$ over $V$ is a
1-manifold with boundary which consists of three 
parts:
(1) the zero set of $\sigma_0$ 
in  $V_{\gamma_1, \cdots, \gamma_k}$; 
 (2) the zero set of $\sigma_1$ 
in  $V'_{\gamma_1, \cdots, \gamma_k}$; 
(3) the zero set of $\sigma$ in the intersection of $V$ with
 $\partial^1 \M^{[m]}_{X_1}(\s_1, a_{(2n+k)}) \times  [0, 1]$, where
$\partial^1 \M^{[m]}_{X_1}(\s_1, a_{(2n+k)})$ is the codimension one boundary
of $\M^{[m]}_{X_1}(\s_1, a_{(2n+k)})$. 

By choosing the base points
$x$ in $X_1$ away from the cylindrical
 end for the $U(1)$-fibration $\tilde \M^{[m]}_{X_1}(\s_1,  a_{(2n+k)})$, 
and noticing that $[x]$ is the generator of $H_0(X_1, \Z)$, 
we know that the contribution of (3) from the codimension one boundary
of $\M^{[m]}_{X_1}(\s_1,  a_{(2n+k)})$ times $[0, 1]$ only comes from
\[
\bigcup_{a_{(2n+k-1)}}\bigcup_{
m_1\in \Z}\M^{[m_1]}_{X_1}(\s_1, a_{(2n+k-1)}) \times 
\hat\M^{[m-m_1]} (a_{(2n+k-1)}, a_{(2n+k)})
\]
 where $a_{(2n+k-1)}$ runs  
over all the possible critical points with  
index $\ind_{X_1}(a_{(2n+k-1)}) =2n+k-1$.  
Now the
transversality condition over the boundary strata over
the intersection of $V$ with $ \M^{[m_1]}_{X_1}(\s_1, a_{(2n+k-1)})  
   \times [0, 1]$ 
provides a set of finitely many oriented points. 
Counting these points gives a number,  denoted
 by $H_{ a_{(2n+k-1)};\  m_1 }$. Then the difference of
the two choices defines a 
number
\[ 
\sum_{a_{(2n+k-1)}}
H_{ a_{(2n+k-1)};\ m_1}\#(\M^{[m-m_1]} (a_{(2n+k-1)}, a_{(2n+k)}) ) ,
\] 
thus,
\[\begin{array}{c}
 \sum_{a_{(2n+k-1)}\atop{m\in \Z; m_1\in \Z} }
H_{ a_{(2n+k-1)};\ m_1}\# (\M^{[m-m_1]} (a_{(2n+k-1)}, a_{(2n+k)}))
t^m a_{(2n+k)} \\[2mm]
  = \partial \bigl( \sum_{a_{(2n+k-1)};m_1\in \Z }
H_{ a_{(2n+k-1)};\ m_1} t^{m_1} a_{(2n+k-1)}\bigr),\end{array}
\]
 is homologous to zero. 
This  
proves that the two choices define the same 
homology class in $HF_{\Gamma, *}^{SW}(Y_0, \t, \Q((t)))$.
\end{proof}

Consider a closed 4-manifold $(X, \s)$  with $b^+_2 \geq 1$ and a $\spinc$ 
structure $\s$. Suppose that  $X$ has 
a decomposition 
$X_1 \cup_{Y_0}  X_2$ 
along a 3-dimensional submanifold $Y_0$  
with $b_1(Y_0)>0$, and with a $\spinc$ structure $\s|_{Y_0}=\t$. We
assume that the class
$c_1(\s|_{Y_0}) = c_1(\t)$ is 
a torsion element in $H^2(Y_0, \Z)$.  
Here we need to  choose a metric  $g_T$ on $X$  
such that $X$ contains a very  long `neck' $ Y_0\times [-T, T]$.   
For $b_2^+ =1$, when the Seiberg--Witten invariant has a chamber  
structure, our choice of the perturbation  compatible
with $\Gamma$ gives rise to 
a fixed chamber in the space of metrics and perturbations. 
Notice that, since we are dealing with the case of a
torsion element $c_1(\t)$, in the case where $b_1(Y_0)=1$ the chamber
also depends on the class $[*\rho]$ in $H^2(Y_0,\R)$ 
Now we can state the gluing formulae of the 
Seiberg--Witten invariants on $X$ as follows. 
 
\begin{The}\label{Glue:trivial} 
Let $X$ be a closed manifold with $b^+ \geq 1$ which is 
written as $X=X_1 \cup_{Y_0} X_2$, with $X_1$ and $X_2$ two $4$-manifolds  
with boundary, 
and with $\partial X_1=-\partial X_2=Y_0$. 
Suppose that we have $\spinc$ structures  
$\s_1$ and $\s_2$, 
on $X_1$ and $X_2$ respectively, such that 
$\s_1|_{Y_0 }\cong \s_2|_{Y_0} \cong \t$ is a torsion $\spinc$ structure 
on $Y_0$. Suppose given  $1$-cycle $\Gamma$ in $ Y_0$
bounds two relative $2$-cycles  
$\Sigma_i$ in $ X_i$ respectively. Then 
the relative Seiberg-Witten invariants $SW_{X_1}(\s_1)$
and $SW_{X_2}(\s_1)$ take values in 
$HF_{\Gamma, *}^{SW} (Y_0, \t, \Q((t)))$ and
$HF_{-\Gamma, *}^{SW} (- Y, -\t, \Q((t)))$ respectively.
Let $\Sigma=\Sigma_1
+\Sigma_2 \in H_2(X, \Z)$, under the natural pairing (\ref{pairing}),  
we have the following  
gluing formula for $z_i \in \A(X_i)$, $i=1,2$, 
\[ 
    \sum_{ \{\s: \s|_{X_i}=\s_i \} }\hspace{-5mm} SW_{X,\s} 
   (z_1z_2)t^{c_1(\s)\cdot \Sigma} =  
   \la SW_{X_1}(\s_1,z_1), SW_{X_2}(\s_2,z_2)\ra. 
\] 
When $b^+=1$, the Seiberg--Witten invariants correspond
to a fixed chamber determined by the perturbation which is compatible
with $\Gamma$. 
\end{The} 
\begin{proof} 
Let $i_k$ denote 
the boundary embedding maps of $Y_0$ into $X_k$. Then 
the set of $\spinc$-structures on $X$
 which agree with $\s_i$ over $X_i$, denoted by $\spinc(X; \s_1,
\s_2)$, form  an affine space over  
 \[ 
H^1(Y_0,  \Z ) /\bigl(Im i_1^* + Im i_2^*\bigr), 
\] 
with $i_1^*: H^1(X_1, \Z) \to H^1(Y_0, \Z)$ and $i_2^*: 
H^1(X_2, \Z) \to H^1(Y_0, \Z)$. 
For any $\s \in \spinc(X; \s_1, \s_2)$, after a generic perturbation, 
the moduli space $\M_X(\s)$ is an oriented, compact, 
 smooth manifold with dimension given by  
\[ 
\frac 14 (c_1(\s)^2 - (2 \chi + 3\sigma )) = d_X(\s) \ge 0. 
\] 
The Seiberg--Witten invariant is a linear functional 
\[ 
SW_X(\s, \cdot): \qquad \AAA(X) \longrightarrow \Z  
\] 
where $\AAA(X) = Sym^*(H_0(X))\otimes \Lambda^*(H_1(X))$.  We summarize 
this definition briefly. We shall only consider monomials
$z=[x]^n \gamma_1 \wedge \cdots \wedge \gamma_k \in \AAA(X)$  
with $2n+k =d_X (\s)$. For any other degree of $z$, we assign
$SW_X(\s, z)$ to be zero. For any monomial as above, 
we can choose smooth 1-dimensional submanifolds representing 
$\gamma_1, \cdots, \gamma_k$.   
Then the holonomy along these loops 
defines a map:  $\M_X(\s) \mapsto U(1)^k$, whose 
generic fiber is a  closed codimension  
 $k$ submanifold $V_{\gamma_1, \cdots, \gamma_k}$ 
in $\M_X(\s)$.   The based monopoles 
define a $U(1)$-fiber bundle $\tilde \M_X (\s) $ over $\M_X(\s)$. 
Then  $SW_X(\s, z)$ is defined to be the  
result of integration of the $n^{th}$ power of  
the first Chern class of this $U(1)$ bundle 
over $V_{\gamma_1, \cdots, \gamma_k}$. Equivalently, 
we can consider the associated rank $n$ complex vector bundle $ E$ 
over the 2n-dimensional manifold  
$V_{\gamma_1, \cdots, \gamma_k}$.  
Then $SW_X(\s, z)$ is obtained  
by counting points with the orientation 
in the zero set of a generic strata transverse section of $E$ over  
$V_{\gamma_1, \cdots, \gamma_k}$.

We can adapt the gluing theory for Seiberg--Witten monopoles developed
in \cite{MW} to the setting of \cite{Fuk} and \cite{Taubes}.
Thus, when $T$ is sufficiently large, 
we obtain a gluing theorem for Seiberg--Witten 
monopoles, which  identifies the 4-dimensional monopoles on $X$  
with the following product: 
\ba 
\bigcup_{\s \in \spinc (X, \s_1, \s_2)}\M_X(\s) \cong  \bigcup_{ a}  
\M_{X_1}(\s_1,  a ) \times \M_{X_2}(\s_2,  a ). 
\label{glue:4d} 
\na 
Notice that, for $[A, \Psi] \in \M_X(\s)$, we  
write $[A, \Psi]$ as $[A_1, \Psi_1] \#^T [A_2, \Psi_2]$ 
under the gluing map (\ref{glue:4d}). We have
\[ 
\E_{\Sigma_1} ([A_1, \Psi_1] )  
+ \E_{\Sigma_2} ([A_2, \Psi_2] )   
= c_1(\s)\cdot \Sigma. 
\] 
Now the gluing formulae follow from the definition 
of the relative invariants $SW_{X_1}(\s_1, z_1)$ and $SW_{X_2}(\s_2,
z_2)$ for   
$z_i \in \AAA (X_i)$, with   
values in $HF^{SW}_{\Gamma, *} (Y_0, \t, \Q((t)))$ and 
 $HF^{SW}_{-\Gamma, *}(- Y_0, - \t)$ respectively. 

Note that, in the case of $b_2^+ (X) =1$, the choice of the 1-cycle
$\Gamma$ fixes the choice of the chamber. In fact, the moduli space
$\M_{Y_0}(\t, \Gamma)$, with the perturbation $[*\rho ] =\eta
PD(\Gamma)$,
contains no reducibles. Thus, in particular, since the relative
invariants take values in $HF_{\Gamma, *}^{SW}(Y_0, \t, \Q((t)))$, 
once $\Gamma$ is fixed, 
there are no wall crossing terms for $SW_{X, \s}$ when we change the 
perturbation and the metric.
\end{proof} 
 
We conclude with the following example.

\begin{Exa}(\cite{vicente}) Let $(Y_0, \s_0)$ be the 3-manifold
$\Sigma_g\times S^1$  
with the trivial $\spinc$ structure $\s_0$, where $\Sigma_g$ is a 
closed Riemann surface of genus $g$.  
Then the gluing formulae in Theorem \ref{Glue:trivial} 
for 4-manifold $\Sigma_g \times S^2$  
can be applied to study  
$HF^{SW}_{\Gamma, *} (\Sigma_g\times S^1, \s_0; \Q((t)))$ and its ring
structure.  
Here we have $\Gamma= [pt \times S^1]$, with the surface $\Sigma$ of
Theorem \ref{Glue:trivial} corresponding to $pt \times S^2$.
It turns out that, as a $Sp(2g, \Z)$-equivariant vector 
space, we have
\[ 
HF^{SW}_{\Gamma, *} (\Sigma_g\times S^1, \s_0; \Q((t))) 
\cong H_*(Sym^{g-1}(\Sigma), \Q) \otimes \Q((t)). 
\] 
In particular, when $\Sigma_g$ is a torus, then 
\[ 
HF^{SW}_{\Gamma, *} (T^2 \times S^1, \s_0; \Q((t))) 
\cong \Q((t)). 
\] 
Upon choosing  the metric on $T^2\times S^2$ such that 
it has a long neck around $T^2 \times S^1$,
 the Seiberg-Witten invariant is computed as the wall-crossing term
from the chamber corresponding to the unperturbed equations
for the positive scalar curvature  which has trivial Seiberg-Witten
invariant, and the other chamber determined by the choice of
$\Gamma$, then we obtain the following result from the
gluing formulae:   
\[\begin{array}{lll} 
&&\la  
 SW_{T^2 \times D^2 }(\s_0, 1), SW_{T^2\times D^2}(\s_0, 1)\ra\\[2mm] 
&=&\sum_{n\geq 1} 
SW_{T^2\times S^2, 2nPD([T^2])} (1) t^{ 2n}\\[2mm] 
&=& \sum_{n\geq 1} n t^{ 2n}. 
\end{array} 
\] 
Thus, up to sign, we can express the relative 
Seiberg--Witten invariant of $T^2\times D^2$ with the trivial $\spinc$ 
structure as 
\[ 
 SW_{T^2 \times D^2 }(\s_0, 1) 
= \displaystyle{\frac 
{1}{t - t^{-1}}}. 
\] 
For $z \in \AAA (T^2 \times D^2)$ of non-zero degree, we have 
$SW_{T^2 \times D^2 }(\s_0, z) =0$. 
\end{Exa}

\vskip .3in 
 
\noindent {\bf Matilde Marcolli}, Department of Mathematics, 
Massachusetts Institute of Technology, 2-275. Cambridge, MA 02139, USA 
\par
\noindent Max-Planck-Institut f\"ur Mathematik, 
Vivatsgasse 7, D-53111 Bonn, 
Germany.  \par  
\noindent matilde\@@math.mit.edu, marcolli\@@mpim-bonn.mpg.de  
 
\vskip .2in 
 
\noindent {\bf Bai-Ling Wang}, Department of Pure Mathematics, 
University of Adelaide, Adelaide SA 5005, Australia. \par 
\noindent bwang\@@maths.adelaide.edu.au 
\par 
\noindent Max-Planck-Institut f\"ur Mathematik, 
Vivatsgasse 7, D-53111 Bonn, 
Germany.  \par 
\noindent bwang\@@mpim-bonn.mpg.de

\end{document}